# SYSTEMATIC SCAN FOR SAMPLING COLORINGS[1]

By Martin Dyer, Leslie Ann Goldberg and Mark Jerrum

*University of Leeds, University of Warwick and University of Edinburgh*

We address the problem of sampling colorings of a graph $G$ by Markov chain simulation. For most of the article we restrict attention to proper $q$-colorings of a path on $n$ vertices (in statistical physics terms, the one-dimensional $q$-state Potts model at zero temperature), though in later sections we widen our scope to general "$H$-colorings" of arbitrary graphs $G$. Existing theoretical analyses of the mixing time of such simulations relate mainly to a dynamics in which a random vertex is selected for updating at each step. However, experimental work is often carried out using systematic strategies that cycle through coordinates in a deterministic manner, a dynamics sometimes known as *systematic scan*. The mixing time of systematic scan seems more difficult to analyze than that of random updates, and little is currently known. In this article we go some way toward correcting this imbalance. By adapting a variety of techniques, we derive upper and lower bounds (often tight) on the mixing time of systematic scan. An unusual feature of systematic scan as far as the analysis is concerned is that it fails to be time reversible.

**1. Introduction.** Many models in statistical physics come under the heading of "spin systems." Such a system is specified by a graph $G$, in our case finite. Configurations of the system are assignments of "spins" to the vertices of $G$. There are assumed to be $q$ possible spins, and, hence, potentially $q^n$ configurations, where $n$ is the number of vertices of $G$, though some of these configurations may be illegal. Each configuration has an energy that comes from summing, over all edges of $G$, the interaction energies between adjacent spins. These energies specify a probability distribution, called the *Boltzmann*

Received June 2004; revised April 2005.

[1]Supported in part by the IST Programme of the EU under contract number IST-1999-14036 (RAND-APX), and by the UK EPSRC through the research grant *Discontinuous Behaviour in the Complexity of Randomized Algorithms*.

*distribution*, on configurations. The Potts model and the hard-core lattice gas model are examples of spin systems.

In this paper for consistency with previous literature, we shall refer to spins as *colors* and to configurations as *states*. Sampling from the Boltzmann distribution is a challenging computational task. Often, the only feasible way of going about it is to simulate a suitable random "dynamics" on configurations. The dynamics has the property of converging to a stationary distribution which is the Boltzmann distribution. This is usually straightforward to arrange. The hard part is proving that the dynamics is "rapidly mixing," that is, converges rapidly to stationarity.

Identifying the vertices of $G$ with the integers $\{1, 2, \ldots, n\}$, we may think of the state space as having coordinates. There is a substantial body of literature concerned with bounding *mixing time* (i.e., time to convergence to near-stationarity) of systems such as those described above. Almost all this theoretical work relates to random single-site updates, which choose a random coordinate for updating at each transition. We shall refer to this strategy as *Glauber dynamics*. (The term "Glauber dynamics" appears not to have a precise agreed meaning. Here we are using the term to signify *single site* updates performed in a *random sequence*. These are certainly aspects of the dynamics first considered by Glauber [18].) However, experimental work is often carried out using systematic strategies that cycle through coordinates in a deterministic manner, a dynamics we refer to as *systematic scan* (or just "scan" for short). The mixing time of systematic scan seems more difficult to analyze that that of Glauber, and little is currently known.

In this paper we take some first steps in analyzing systematic scan for spin systems. Our setting will be very simple; indeed, for the most part, we will restrict attention to proper $q$-colorings of a path of $n$ vertices (in statistical physics terms, the one-dimensional $q$-state Potts model at zero temperature). To compensate for the simple setting, we provide tight (i.e., matching within a constant factor) upper and lower bounds on mixing time. Measuring mixing time in terms of the number of updates of individual vertices (so that one scan equates to $n$ updates), we show that when $q = 3$, mixing occurs in $\Theta(n^3 \log n)$ updates, whether Glauber dynamics or systematic scan is used; while when $q \geq 4$, mixing occurs in $\Theta(n \log n)$ updates, again independently of whether Glauber or scan is used. Our main tools are harmonic analysis [29], path coupling [6] and disagreement percolation [28].

Later in the paper we considerably widen the setting from usual proper colorings to general $H$-colorings (also known a graph homomorphisms), but staying at first with the path as the underlying graph. $H$-colorings model arbitrary spin systems with symmetric "hard" constraints. We show that, for any $H$, Glauber mixes in $O(n^5)$ updates and scan in $O(n^6)$ updates. The former bound is unlikely to be tight, and the latter even less so. The method here is that of canonical paths [10, 27].

Finally, we consider $H$-colorings of a general graph $G$, and compare the mixing times of scan and Glauber. We show that, for any $H$, these are within a polynomial factor of each other (in terms of total number of individual updates performed), at least when $G$ is of bounded degree. The question of whether scan can ever be faster than Glauber, or vice versa, remains a tantalizing open problem. The only situation where a gap is known is the rather uninteresting one that arises when $G$ is the empty graph, where Glauber requires $\Theta(n \log n)$ updates [13], while scan clearly mixes in one sweep.

1.1. *Previous work.* Amit [3] has investigated systematic scan in the context of sampling from multivariate Gaussian distributions. In this instance, one iteration of systematic scan applies a "heat-bath" update to each coordinate axis in turn. Amit precisely calculates the spectral gap of the scan operator and, hence, bounds the mixing time. He also estimates the spectral gap of a similar process on perturbed Gaussian distributions.

In another application of systematic scan—this time more combinatorial in nature and slightly closer to the one studied here—Diaconis and Ram [8] consider the problem of generating random elements of a finite group. The *systematic scan* Metropolis algorithm cycles through the generators in order, and flips coins to decide whether or not to multiply by each generator in turn. The *random update* algorithm chooses one of the $n$ generators uniformly at random at each step. For the symmetric group, they show that the systematic scan algorithm mixes in $\Theta(n)$ scans, so consideration of $\Theta(n^2)$ selections of generators is necessary and sufficient for mixing. They consider two different scanning strategies from [17]—the same results hold for both strategies. Matching results (in terms of the number of generators considered) are given by Benjamini et al. [5] for the random update strategy. Diaconis and Ram also consider the hypercube and the dihedral group. For the hypercube, they show that $\Theta(n \log n)$ updates are necessary and sufficient, whether one is doing random updates or systematic scan. For the dihedral group, both strategies take $\Theta(n)$ updates. Diaconis and Ram point out that careful analysis of rates of convergence for the Metropolis algorithm is completely open in nongroup cases.

For a brief review of other work on systematic scan, consult Diaconis and Ram [8], Section 2b.

**2. Definitions and notation.** The *variation distance* between distributions $\theta_1$ and $\theta_2$ on $\Omega$ is

$$\mathrm{d_{TV}}(\theta_1, \theta_2) = \tfrac{1}{2} \sum_i |\theta_1(i) - \theta_2(i)| = \max_{A \subseteq \Omega} |\theta_1(A) - \theta_2(A)|.$$

For a discrete ergodic Markov chain $\mathcal{M}$ with transition matrix $P$ and stationary distribution $\pi$, and a specified initial state $x$, the mixing time (as a function of the deviation $\varepsilon$ from stationarity) is

$$\mathrm{Mix}_x(\mathcal{M},\varepsilon) = \min\{t > 0 : \mathrm{d_{TV}}(P^t(x,\cdot),\pi(\cdot)) \le \varepsilon\}.$$

The mixing time of $\mathcal{M}$ is $\mathrm{Mix}(\mathcal{M},\varepsilon) = \max_x \mathrm{Mix}_x(\mathcal{M},\varepsilon)$.

Suppose $G$ is an undirected graph with vertex set $\{1,\dots,n\}$. To avoid trivialities, we assume $n > 3$. We consider $q$-colorings of $G$, where $q \ge 3$. Formally, a coloring $\sigma$ is a vector $\sigma = (\sigma_1,\dots,\sigma_n)$ in which $\sigma_i \in \{0,\dots,q-1\}$ denotes the color of vertex $i$. A coloring is *proper* if adjacent vertices receive different colors. $\Omega^+ = \{0,\dots,q-1\}^n$ is the set of all colorings (proper and improper), while $\Omega$ is the set of all proper colorings.

A Markov chain with state space $\Omega$ starts at a coloring $\sigma(0)$ and visits a sequence of colorings $\sigma(0),\sigma(1),\dots$. We often use $\tau$ to denote a coloring (when we need two names). The two Markov chains that we study are as follows:

- $\mathcal{M}_{\mathrm{Gl}}$ (Glauber): Choose vertex $v$ uniformly at random; do Metropolis($v$).
- $\mathcal{M}_{\rightarrow}$ (Systematic scan): For $v := 1$ to $n$, do Metropolis($v$).

The procedure Metropolis($v$) used in both of the above dynamics performs as follows: A color $c$ is chosen uniformly at random. A proposed new coloring is formed by recoloring vertex $v$ with color $c$. This proposed move is accepted if and only if color $c$ is not used at any neighbor of $v$.

Let $P_{\mathrm{Gl}}$ be the transition matrix of $\mathcal{M}_{\mathrm{Gl}}$ and Let $P_{\rightarrow}$ be the transition matrix of $\mathcal{M}_{\rightarrow}$. It will be convenient in our proofs to consider reverse systematic scan:

- $\mathcal{M}_{\leftarrow}$ (Reverse scan): For $v := n$ down to 1, do Metropolis($v$).

Let $P_{\leftarrow}$ be the transition matrix of $\mathcal{M}_{\leftarrow}$. Observe that $\mathcal{M}_{\leftarrow}$ is the time reversal of $\mathcal{M}_{\rightarrow}$, since $P_{\rightarrow}(\sigma,\sigma') = P_{\leftarrow}(\sigma',\sigma)$ for all $\sigma,\sigma' \in \Omega$.

Let $\mathcal{M}$ be any discrete Markov chain with transition matrix $P$, stationary distribution $\pi$ and state space $\Omega$. Define the *optimal Poincaré constant* of $\mathcal{M}$ by

$$\lambda(\mathcal{M}) = \inf_{f:\Omega\to\mathbb{R}} \frac{\mathcal{E}_{\mathcal{M}}(f,f)}{\mathrm{var}_\pi(f)},$$

where the inf is over all nonconstant functions from $\Omega$ to $\mathbb{R}$ and the *Dirichlet form* is given by

$$\mathcal{E}_{\mathcal{M}}(f,f) = \tfrac{1}{2} \sum_{x,y\in\Omega} \pi(x)P(x,y)(f(x)-f(y))^2$$

and

$$\mathrm{var}_\pi(f) = \sum_{x\in\Omega} \pi(x)(f(x)-\mathbf{E}_\pi f)^2 = \tfrac{1}{2} \sum_{x,y\in\Omega} \pi(x)\pi(y)(f(x)-f(y))^2.$$

If $\mathcal{M}$ is time-reversible with respect to $\pi$ [i.e., $\pi(\sigma)P(\sigma,\tau) = \pi(\tau)P(\tau,\sigma)$], then the eigenvalues of $P$ are real and can be written $1 = \beta_0 \geq \beta_1 \geq \cdots \geq \beta_{|\Omega|-1} \geq -1$. Then $\lambda(\mathcal{M})$ is is equal to $1 - \beta_1$.

Some of our rapid-mixing proofs will use the method of *path coupling* [6]. In our path-coupling proofs, we will define partial couplings on the set $S$, which will always be the set of pairs of colorings that differ on a single vertex.

For most of the paper we consider the case in which $G$ is a path going left to right from vertex 1 to vertex $n$. Kenyon and Randall [24] have shown that, for every $q$, the *block dynamics*, which updates a sufficiently large constant-length path at each step, mixes in time $O(n \log n)$. Our results show that this upper bound holds for single-site dynamics for $q \geq 4$, but not for $q = 3$.

In our analysis for $q = 3$ we will study two auxiliary Markov chains on state space $\Upsilon = \{-1, 1\}^{n-1}$. A configuration $X \in \Upsilon$ is a vector $X = (X_1, \ldots, X_{n-1})$. The corresponding Markov chain evolves as $X(0), X(1), \ldots.$ The next section generalizes this framework.

**3. The analysis technique for $q = 3$.** The following develops an idea of Wilson [29] for lower bounding the convergence rate of certain types of Markov chains.

Let $\mathcal{M}$ be a finite ergodic Markov chain with transition matrix $P$ and state space $\Upsilon \subseteq \mathbb{Z}^m$. (This is a more general setting than our current application demands, but it is the natural one in which to develop the ideas.) Suppose there exists a matrix $A$ such that $\mathbf{E}[X(1)|X(0)] = AX(0)$ for all $X(0) \in \Upsilon$. (The method may still be applicable when we have only an affine dependence here. For provided $A - I$ is invertible, an affine dependence $\mathbf{E}[X(1)] = AX(0) + b$ can be reduced to one of the required form by moving the origin in $\Upsilon$. In particular, $\mathbf{E}[X(1)] = AX(0) + b$ is the same as $\mathbf{E}[X(1) + c] = A(X(0) + c)$ for $b = (A - I)c$.) We will assume that $A$ has real eigenvalues, though it is possible to extend the method to complex eigenvalues. We may further assume that $A$ has only nonnegative eigenvalues, since otherwise we can consider the two-step chain $\mathcal{M}^2 = (\Upsilon, P^2)$ which converges exactly twice as fast. Now let $\lambda$ be any eigenvalue of $A$, with left eigenvector $w$. Then,

$$\mathbf{E}[wX(t)|X(0)] = wA^t X(0) = \lambda^t wX(0). \tag{1}$$

Let $\Phi_t = wX(t)$. To obtain the strongest lower bound, we choose $\lambda$ to be the largest eigenvalue such that there exist $x, y \in \Upsilon$ with $wx \neq wy$. Then we choose $X(0) = \arg\max_x |wx|$. Since $w$ is defined only to scalar multiples, we may assume $wX(0) > 0$. It follows from (1) that $\lambda \leq 1$. Otherwise $\limsup_{t\to\infty} \mathbf{E}[\Phi_t] = \infty$, contradicting the finiteness of $\mathcal{M}$. If $\lambda = 1$, we have $\mathbf{E}[\Phi_t] = \Phi_0$ for all $t$. But $\Phi_t \leq \Phi_0$, so we must have $\Phi_t = \Phi_0$ for all $t$. Using ergodicity of $\mathcal{M}$, this implies $wx = wy$ for all $x, y \in \Upsilon$, contradicting our choice of $\lambda$. Thus, $\lambda < 1$ and, hence, $\lim_{t\to\infty} \mathbf{E}[wX(t)] = 0$. If $X(\infty)$ denotes (a r.v. with) the equilibrium distribution, it follows that $\mathbf{E}[wX(\infty)] = 0$.

We will now consider the quantities $\mathbf{E}[\Phi_t|\Phi_{t-1}]$ and $\operatorname{var}(\Phi_t|\Phi_{t-1})$. Definitions of conditional expectations and variances can be found in [11] (pages 190–198). We will use the fact that $\operatorname{var}(Y) = \mathbf{E}[\operatorname{var}(Y|X)] + \operatorname{var}(\mathbf{E}[Y|X])$ (page 198). Suppose that $\mathbf{E}[\operatorname{var}(\Phi_t|\Phi_{t-1})] \le \rho$ for all $t > 0$, and let $\nu = \rho/(1-\lambda^2)$. Now using $\mathbf{E}[\Phi_t|\Phi_{t-1}] = \lambda\Phi_{t-1}$ and $\operatorname{var}(\Phi_0) = 0$,

$$\begin{aligned} \operatorname{var}(\Phi_t) &= \mathbf{E}[\operatorname{var}(\Phi_t|\Phi_{t-1})] + \operatorname{var}(\mathbf{E}[\Phi_t|\Phi_{t-1}]) \\ &= \mathbf{E}[\operatorname{var}(\Phi_t|\Phi_{t-1})] + \operatorname{var}(\lambda\Phi_{t-1}) \\ &= \mathbf{E}[\operatorname{var}(\Phi_t|\Phi_{t-1})] + \lambda^2 \operatorname{var}(\Phi_{t-1}) \qquad (2)\\ &\le \rho + \lambda^2 \operatorname{var}(\Phi_{t-1}) \\ &\le \sum_{i=0}^{t-1} \lambda^{2i}\rho < \rho/(1-\lambda^2) = \nu. \end{aligned}$$

Instead of (2), Wilson uses $\nu = R/2\gamma$, where $\gamma = 1-\lambda$ and $\mathbf{E}[(\Phi_t - \Phi_{t-1})^2|\Phi_{t-1}] \le R$. The calculation to justify this is longer, and the conclusion is not valid for all $\lambda$. However, since $\rho \le R$ and usually $\lambda = o(1)$, (2) implies Wilson's bound asymptotically, but, in general, they are incomparable. Now, using Chebyshev's inequality,

$$\Pr\left(\Phi_t < \lambda^t\Phi_0 - \sqrt{\frac{2\nu}{\varepsilon}}\right) < \frac{1}{2}\varepsilon \quad \text{and} \quad \Pr\left(\Phi_\infty > \sqrt{\frac{2\nu}{\varepsilon}}\right) < \frac{1}{2}\varepsilon.$$

Thus, $d_{\mathrm{TV}}(\Phi_t, \Phi_\infty) \le 1-\varepsilon$ only if $\lambda^t\Phi_0 < 2\sqrt{2\nu/\varepsilon}$. [We will abuse the notation $d_{\mathrm{TV}}(\cdot,\cdot)$ for variation distance by extending it to random variables.] The latter inequality holds only if

$$t > \frac{\ln(\sqrt{\varepsilon/8}\,\Phi_0/\sqrt{\nu})}{\ln(1/\lambda)} \ge \frac{\lambda\ln(\sqrt{\varepsilon/8}\Phi_0/\sqrt{\nu})}{1-\lambda}.$$

Setting $\varepsilon = \frac{1}{2}$, we find that

$$\operatorname{Mix}\left(\mathcal{M}, \frac{1}{2}\right) \ge \frac{\lambda\ln(\Phi_0/4\sqrt{\nu})}{1-\lambda}. \qquad (3)$$

We say that a Markov chain is *monotone* with respect to a partial order $\le$ on its state space if two realizations $X(t)$ and $Y(t)$ of it may be coupled so that $X(0) \ge Y(0)$ implies $X(t) \ge Y(t)$ for all $t \in \mathbb{N}$. We refer to such a coupling as a "monotone coupling." Suppose $\mathcal{M}$ is monotone with respect to the product partial order $\le$ on $\mathbb{R}^m$, that is, the partial order defined by $x \le y$ if and only if $x_i \le y_i$ for all $i \in \{1,\dots,m\}$. If the weight vector $w > 0$ (in the product order), we can use it to bound the mixing time from above. Let us re-scale $w$ so that $\min_i w_i = 1$. Let $d(x,y) = \sum_{i=1}^{m} w_i |x_i - y_i|$ for $x, y \in \Upsilon$. Then $d$ is a metric, since $w > 0$. Now consider $x, y \in \Upsilon$ with $x \geq y$ and let $X(t), Y(t)$ be a monotone coupling with $X(0) = x$ and $Y(0) = y$. Since $X(t) \geq Y(t)$,

$$\begin{aligned}\mathbf{E}[d(X(t),Y(t))] &= \mathbf{E}[w(X(t)-Y(t))] = wA(X(t-1)-Y(t-1)) \\ &= \lambda w(X(t-1)-Y(t-1)) = \lambda\, d(X(t-1),Y(t-1)).\end{aligned}$$

So

$$\begin{aligned}\mathrm{d}_{\mathrm{TV}}(X(t),Y(t)) &\leq \Pr[X(t) \neq Y(t)] \\ &\leq \mathbf{E}[d(X(t),Y(t))] \\ &\leq \lambda^t d(X(0),Y(0)) \leq 2\lambda^t \Phi_0,\end{aligned} \tag{4}$$

where the final step is by the triangle inequality. Thus, $\mathrm{d}_{\mathrm{TV}}(X(t),Y(t)) \leq \varepsilon$ holds, provided $t \geq \ln(2\Phi_0/\varepsilon)/\ln(1/\lambda)$, that is, provided $t \geq \ln(2\Phi_0/\varepsilon)/(1-\lambda)$.

We would like to draw a similar conclusion when $x$ and $y$ are incomparable. We can do so provided the state space contains states $\top$ and $\bot$ satisfying $\bot \leq z \leq \top$ for all $z \in \Upsilon$. In this case, $\Pr(X(t) \neq Y(t) | X(0) = x, Y(0) = y) \leq \Pr(X(t) \neq Y(t) | X(0) = \top, Y(0) = \bot)$ so we can apply (4) starting from $X(0) = \top$ and $Y(0) = \bot$. Thus,

$$\mathrm{Mix}(\mathcal{M}, \varepsilon) \leq \ln(2\Phi_0/\varepsilon)/(1-\lambda). \tag{5}$$

When $\nu$ is sufficiently small with respect to $\Phi_0$, the upper bound (5) and the lower bound (3) agree to within a constant factor on the time to reach variation distance $\frac{1}{2}$, say.

3.1. *Bounding* $\mathbf{E}[\mathrm{var}(\Phi_t|\Phi_{t-1})]$. In order to use the technique in Section 3, we have to find a $\rho$ such that $\mathbf{E}[\mathrm{var}(\Phi_t|\Phi_{t-1})] \leq \rho$, where the expectation is over $\Phi_{t-1}$.

For this, let $Z_t = \Phi_t - \Phi_{t-1}$. Then

$$\begin{aligned}\mathbf{E}[\mathrm{var}(\Phi_t|\Phi_{t-1})] &= \mathbf{E}[\mathrm{var}(\Phi_{t-1} + Z_t|\Phi_{t-1})] \\ &= \mathbf{E}[\mathrm{var}(Z_t|\Phi_{t-1})] \\ &= \mathbf{E}[\mathbf{E}[Z_t^2|\Phi_{t-1}] - (\mathbf{E}[Z_t|\Phi_{t-1}])^2] \\ &\leq \mathbf{E}[\mathbf{E}[(\Phi_t - \Phi_{t-1})^2|\Phi_{t-1}]] \\ &\leq \max_{\Phi_{t-1}} \mathbf{E}[(\Phi_t - \Phi_{t-1})^2|\Phi_{t-1}].\end{aligned}$$

We will use the above inequality to find a suitable $\rho$.

3.2. *A benchmark example.* We illustrate the technique by applying it to a simple example whose analysis was also given in the Introduction to [8]. Consider mixing on the cube $\{-1,+1\}^m$ of the chain which changes the sign of a uniform random coordinate with probability $\frac{1}{2}$. Then

$$\mathbf{E}[X_i(t+1)] = (1-1/m)X_i(t),$$

so $A=(1-1/m)I$, and all its eigenvalues are equal to $(1-1/m)$. We may choose an arbitrary $w$, say, the vector $(1,1,\dots,1)$. Then we can take $\rho=2$, so $\nu=2m/(2-1/m)\le 2m$. Taking $X(0)=w$, $\Phi_0=m$, and the lower bound (3) for mixing time is $\frac{1}{2}m\ln m - O(m)$. This chain is monotone, and the upper bound (5) is $m\ln m + O(m)$.

**4. Glauber dynamics for $q=3$ mixes in $\Theta(n^3\log n)$ updates.** Let $G$ be a path going left to right from vertex 1 to vertex $n$. Recall that $\Omega$ is the set of all proper colorings of $G$.

4.1. *Analysis of a related Markov chain.* Let $\sigma$ be a coloring in $\Omega$. Note that, for every $i\in\{1,\dots,n-1\}$, we either have $\sigma_{i+1}=\sigma_i+1 \pmod 3$ or $\sigma_{i+1}=\sigma_i-1 \pmod 3$. We can associate $\sigma$ with a vector $X\in\Upsilon=\{-1,1\}^{n-1}$. $X_i$ is 1 if $\sigma_{i+1}=\sigma_i+1 \pmod 3$ and $X_i=-1$ otherwise. (Note that three colorings are mapped to the same configuration $X$—given $\sigma_1$ and $X$, the coloring $\sigma$ can be recovered.)

The Markov chain $\mathcal{M}_{\mathrm{Gl}}$ can be associated with a Markov chain $\mathcal{M}^{\pm}_{\mathrm{Gl}}$ on $\Upsilon$. The moves of $\mathcal{M}^{\pm}_{\mathrm{Gl}}$ (from configuration $X$) are as follows. Choose $r\in\{1,\dots,n\}$ uniformly at random. If $r=1$ (resp. $r=n$), then either, with probability $\frac{1}{3}$, change the sign of $X_1$ (resp. $X_{n-1}$) or, with the complementary probability, do nothing. Otherwise (i.e., if $1<r<n$) then either, with probability $\frac{1}{3}$, exchange $X_{r-1}$ and $X_r$ or, with the complementary probability, do nothing.

In this section we analyze the mixing rate of $\mathcal{M}^{\pm}_{\mathrm{Gl}}$. Note that this chain is monotone with respect to the usual partial order on $\Upsilon$. Then straightforward calculations give

$$\mathbf{E}[X(t+1)] = AX(t) \qquad \text{where } A = I - \frac{1}{3n}B$$

and

$$B = \begin{bmatrix} 3 & -1 & & & & \\ -1 & 2 & -1 & & & \\ & -1 & 2 & -1 & & \\ & & & \ddots & & \\ & & & -1 & 2 & -1 \\ & & & & -1 & 3 \end{bmatrix}.$$

Note that $A$ is symmetric so has all real eigenvalues. Moreover, $A$ is nonnegative and irreducible, so its largest eigenvalue $\lambda$ has a positive eigenvector $w$. The eigenvectors are identical to those of $B$, and $\lambda = (1 - \lambda'/3n)$, where $\lambda'$ is the smallest eigenvalue of $B$. The "generic" row gives the equation

$$-w_{i-1} + 2w_i - w_{i+1} = \lambda' w_i, \tag{6}$$

the form of which suggests a simple harmonic oscillation. So we will try the solution $w_i = c_n \sin(\alpha i + \beta)$, where $c_n$ is a positive scaling factor to be chosen later. Substituting in (6) gives $\lambda' = 2(1 - \cos\alpha) = 4\sin^2(\alpha/2)$. We also have the two "boundary conditions"

$$3w_1 - w_2 = \lambda' w_1, \qquad -w_{n-2} + 3w_{n-1} = \lambda' w_{n-1}. \tag{7}$$

The first equation in (7) gives $\sin(\alpha + \beta) = -\sin\beta$, that is, $\beta = -\alpha/2$. The second then gives $\sin((n - \frac{1}{2})\alpha) = -\sin((n - \frac{3}{2})\alpha)$, so $(n - \frac{1}{2})\alpha = 2\pi - (n - \frac{3}{2})\alpha$, that is, $\alpha = \pi/(n-1)$. Thus,

$$w_i = c_n \sin\left(\frac{\pi(i - 1/2)}{n-1}\right) > 0, \qquad i = 1, \ldots, n-1,$$

and $w = (w_i)$ is the (positive) eigenvector corresponding to the largest eigenvalue. Our upper bound on mixing time requires $w_i \geq 1$, for all $1 \leq i \leq n-1$, and we set $c_n \sim 2n/\pi$ to achieve this. (The symbol "$\sim$" denotes asymptotic convergence as $n \to \infty$.)

Now, if we let $w_0 = w_n = 0$, we may take

$$\rho = 2 \max_{1 \leq i \leq n} (w_i - w_{i-1})^2 = 2(w_2 - w_1)^2 \sim 8.$$

Also, $\lambda = 1 - 4\sin^2(\pi/(2n-2))/3n$, so $1 - \lambda \sim \pi^2/3n^3$. Hence, $\nu \sim 12n^3/\pi^2$. Taking $X(0)$ to be the all 1's vector,

$$\begin{aligned}\Phi_0 &= c_n \sum_{i=1}^{n-1} \sin\left(\frac{\pi(i-1/2)}{n-1}\right)\\ &= c_n \mathrm{Im}\left[\sum_{j=1}^{n-1} \exp\left(\frac{\mathbf{i}\,\pi(j-1/2)}{n-1}\right)\right]\\ &= c_n \mathrm{cosec}\left(\frac{\pi}{2(n-1)}\right) \sim \left(\frac{2n}{\pi}\right)^2,\end{aligned}$$

where $\mathbf{i} = \sqrt{-1}$ in the second equality and the final equality follows from simplifying the geometric series as follows. For $\ell = i\pi/(n-1)$, the sum is equal to

$$\begin{aligned}e^{-\ell/2}\left(\frac{e^{\ell} - e^{\ell n}}{1 - e^{\ell}}\right) &= \frac{-2}{\exp(i\pi/2(n-1)) - \exp(-i\pi/2(n-1))}\\ &= \frac{i}{\sin(\pi/2(n-1))}.\end{aligned}$$

Substituting for $\lambda$, $\Phi_0$ and $\nu$ in the mixing time lower bound, (3) yields $\mathrm{Mix}(\mathcal{M}_{\mathrm{Gl}}^{\pm}, \frac{1}{2}) \geq \frac{3}{2}\pi^{-2}n^3 \ln n - O(n^3)$. Also (for any positive $\varepsilon$), the upper bound (5) is $\mathrm{Mix}(\mathcal{M}_{\mathrm{Gl}}^{\pm}, \varepsilon) \leq 3\pi^{-2}n^3(2\ln n + \ln \varepsilon^{-1}) + O(n^3)$. In summary, the mixing time of $\mathcal{M}_{\mathrm{Gl}}^{\pm}$ is $\Theta(n^3 \log n)$.

4.2. *Distance measures and a lower bound for Glauber dynamics.* We will use two distance measures to analyze the Glauber-dynamics Markov chain $\mathcal{M}_{\mathrm{Gl}}$. First, we define the distance $d_1(\sigma,\tau)$ for $\sigma \in \Omega$ and $\tau \in \Omega$ follows. Let $X$ be the member of $\Upsilon$ associated with $\sigma$ and $Y$ be the member of $\Upsilon$ associated with $\tau$. Let $d_1(\sigma,\tau) = \mathrm{Ham}(X,Y)$, where $\mathrm{Ham}(X,Y)$ is the Hamming distance between $X$ and $Y$, which is the number of indices $i$ such that $X_i \neq Y_i$.

Using distance measure $d_1$, the lower bound from Section 4.1 applies directly to $\mathcal{M}_{\mathrm{Gl}}$. Thus, we obtain the following theorem.

THEOREM 1. *Let $G$ be the $n$-vertex path, and let $q = 3$. Then a lower bound on the mixing time of the Markov chain $\mathcal{M}_{\mathrm{Gl}}$ on the state space $\Omega$ is given by* $\mathrm{Mix}(\mathcal{M}_{\mathrm{Gl}}, \frac{1}{2}) \geq \frac{3}{2}\pi^{-2}n^3 \ln n + O(n^3)$.

In order to upper-bound the mixing time of $\mathcal{M}_{\mathrm{Gl}}$, we will also define a second distance measure.

Give the vertices $1,\dots,n$ weights $\lambda_1,\dots,\lambda_n$, respectively. These weights are positive rationals. Denote by $S \subset \Omega \times \Omega$ the set of all pairs of states (colorings) that differ at a single vertex (i.e., are Hamming distance 1 apart). If $(\sigma,\tau) \in S$ differs at vertex $i$, then let $\phi(\sigma,\tau) = \lambda_i$. Define the function $d_2$ on $\Omega \times \Omega$ as follows. For each pair $(\sigma,\tau) \in \Omega \times \Omega$, let

$$d_2(\sigma,\tau) = \min_{\omega(0),\dots,\omega(k)} \sum_{j=0}^{k-1} \phi(\omega(j),\omega(j+1)), \tag{8}$$

where the minimum is over all paths $\sigma = \omega(0),\dots,\omega(k) = \tau$ such that each $\omega_j \in \Omega$ and each pair $(\omega(j),\omega(j+1)) \in S$. A path $\omega(0),\dots,\omega(k)$ satisfying (8) is referred to as a *geodesic* path from $\sigma$ to $\tau$.

In our couplings, we will want to be able to bound the expected change in the distance $d_2$. In order to do this, we use height functions. A *height function* $h$ corresponding to a proper coloring $\sigma$ is a vector in $\mathbb{Z}^n$ satisfying the following properties:

1. For every vertex $i$, $h_i \equiv i \pmod 2$.
2. For every vertex $i$, $h_i \equiv \sigma_i \pmod 3$.
3. For every edge $(i,i+1)$, $|h_i - h_{i+1}| = 1$.

The height function is unique up to an additive constant. We define the distance between two height functions, $h$ and $h^*$, to be

$$d(h,h^*) = \sum_{i\in\{1,\ldots,n\}} \frac{|h_i - h_i^*|\,\lambda_i}{2}.$$

Let $\mathcal{H}(\sigma)$ denote the set of height functions corresponding to coloring $\sigma$.

LEMMA 2. *For any pair of colorings* $(\sigma,\tau)\in\Omega\times\Omega$,

$$d_2(\sigma,\tau) = \min_{h\in\mathcal{H}(\sigma),h^*\in\mathcal{H}(\tau)} d(h,h^*).$$

PROOF. To show that

$$d_2(\sigma,\tau) \geq \min_{h\in\mathcal{H}(\sigma),h^*\in\mathcal{H}(\tau)} d(h,h^*),$$

consider a geodesic path from $\sigma$ to $\tau$. Let $h'(0)$ be any height function in $\mathcal{H}(\sigma)$ and let $h'(0),\ldots,h'(k)$ be the sequence of height functions corresponding to the geodesic path. Now

$$\begin{aligned}\min_{h\in\mathcal{H}(\sigma),h^*\in\mathcal{H}(\tau)} d(h,h^*) &\leq d(h'(0),h'(k))\\ &\leq \sum_{i=0}^{k-1} d(h'(i),h'(i+1)) = d_2(\sigma,\tau).\end{aligned}$$

To show that

$$d_2(\sigma,\tau) \leq \min_{h\in\mathcal{H}(\sigma),h^*\in\mathcal{H}(\tau)} d(h,h^*),$$

consider any $h\in\mathcal{H}(\sigma)$ and $h^*\in\mathcal{H}(\tau)$. A "height-function transformation" (see [20]) either takes a local maximum of a height function and pushes it down by two or takes a local minimum and pushes it up by two. We can show that there is a sequence $h = h(0),\ldots,h(k) = h^*$ of height-function transformations transforming $h$ into $h^*$ that chooses each vertex $v$ only $|h_v - h_v^*|/2$ times. (This can be proved by induction on $\sum_v |h_v - h_v^*|$. See Lemma 4.3 of [20].) Now let $\omega(0),\ldots,\omega(k)$ be the sequence of colorings corresponding to $h(0),\ldots,h(k)$. Note that

$$\sum_{j=0}^{k-1} \phi(\omega(j),\omega(j+1)) \leq d(h,h^*).$$

Thus, $d_2(\sigma,\tau) \leq d(h,h^*)$. $\square$

4.3. *An upper bound for Glauber dynamics.* Our upper bound comes from a two-stage argument. In the first stage we observe the evolution of $\mathcal{M}_{\text{Gl}}$, not directly, but via the auxiliary Markov chain $\mathcal{M}^{\pm}_{\text{Gl}}$. We know from Section 4.1 that the latter mixes in time $\Theta(n^3 \log n)$. Each state of $\mathcal{M}^{\pm}_{\text{Gl}}$ corresponds to $q$ states of $\mathcal{M}_{\text{Gl}}$, so at this point we know that $\mathcal{M}_{\text{Gl}}$ has mixed modulo a cyclic permutation of colors. In the second stage we show, using the $d_2$ metric, that two colorings $\sigma$ and $\tau$ differing by such a permutation may be coupled in a further $O(n^3)$ steps. The first stage gives the coupling a head start in the sense that $\sigma$ and $\tau$ are already quite close in the $d_2$ metric. Omitting the first stage and running the $d_2$-coupling in isolation would yield only an $O(n^5)$ bound on mixing time.

Recall that $d_1$ is Hamming distance on $\Upsilon$. Suppose $(\sigma(0), \tau(0)) \in \Omega \times \Omega$. Then

$$\Pr(d_1(\sigma(t), \tau(t)) \geq 1) \leq \mathbf{E}(d_1(\sigma(t), \tau(t))).$$

Applying (4) to the analysis in Section 4.1, the right-hand side is at most $2\lambda^t \Phi_0$, where $\Phi_0 = \Theta(n^2)$ and $1 - \lambda \sim \pi^2/2n^3$. So for some $t' = O(n^3 \log n)$, we will have $d_1(\sigma(t'), \tau(t')) = 0$, with probability at least $\frac{39}{40}$. By Lemma 2, $d_1(\sigma(t'), \tau(t')) = 0$ implies that $d_2(\sigma(t'), \tau(t')) = 0$ or $d_2(\sigma(t'), \tau(t')) = \sum_{i \in \{1,\ldots,n\}} \lambda_i$.

Now choose weights $\lambda_1 = \lambda_n = 1/2$ and $\lambda_2 = \cdots = \lambda_{n-1} = 1$. We use path coupling on pairs $(\sigma(0), \tau(0)) \in S$. Starting with such a pair, run $t'$ steps to get $(\sigma(t'), \tau(t'))$. With probability at least $39/40$, $d_1(\sigma(t'), \tau(t')) = 0$, in which case either $\sigma(t') = \tau(t')$ or $d_2(\sigma(t'), \tau(t')) = n - 1$. If the former holds, we are done, so suppose the latter. We now carry on from $(\sigma(t'), \tau(t'))$ using the identity coupling (i.e., to say the coupling that chooses the same vertex in both copies, and proposes the same color $c$ in both). We will show in Section 4.3.1 below that, if we take any $(\sigma, \tau) \in \Omega \times \Omega$ and produce $(\sigma', \tau')$ by one step of the identity coupling, then $\mathbf{E}[d_2(\sigma', \tau')] \leq d_2(\sigma, \tau)$. Thus, $D_t = d_2(\sigma(t' + t), \tau(t' + t))$ is a super-martingale with $D_0 = n - 1$. In Section 4.3.2 below, we will define a quantity $V = \Omega(1/n)$ and show that, for all $t$ and all values of $D_t$ other than 0, $\mathbf{E}[(D_{t+1} - D_t)^2 | D_t] \geq V$. Let $B = 10n$, and let $T$ be the first time at which either (a) $D_t = 0$ (i.e., coupling occurs), or (b) $D_t \geq B$. Note that $T$ is a stopping time. Define $Z_t = (B - D_t)^2 - Vt$, and observe (see [25]) that $Z_{t \wedge T}$ is a sub-martingale, where $t \wedge T$ denotes the minimum of $t$ and $T$. Let $p$ be the probability that (a) occurs. By the optional stopping theorem $E[D_T] \leq D_0$, so $(1 - p)B \leq E[D_T] \leq D_0$ and $p \geq 1 - D_0/B \geq \frac{9}{10}$. Also, by the optional stopping theorem,

$$\begin{aligned} pB^2 + (1-p)E[(B - D_T)^2 | D_T \geq B] - VE[T] \\ = E[(B - D_T)^2] - VE[T] = E[Z_T] \\ \geq Z_0 = (B - D_0)^2 > 0. \end{aligned}$$

Since $|D_t - D_{t-1}| \le 2$, $(1-p)E[(B-D_T)^2|D_T \ge B] \le 4 < pB^2$ so $\mathbf{E}[T] \le (2pB^2)/V$. Conditioning on (a) occurring, it follows that $\mathbf{E}[T|D_T = 0] \le 2B^2/V$. Hence, $\Pr(T > 20B^2/V|D_T = 0) \le \frac{1}{10}$. So, if we now run the identity coupling for $20B^2/V = O(n^3)$ steps, then $\sigma$ and $\tau$ will fail to couple with probability at most $\frac{1}{40} + 2 \times \frac{1}{10} < \frac{1}{4}$. Thus, we have shown the following.

THEOREM 3. *Let $G$ be the $n$-vertex path, and let $q = 3$. Consider the Markov chain $\mathcal{M}_{\text{Gl}}$ on the state space $\Omega$. Then* $\text{Mix}(\mathcal{M}_{\text{Gl}}, \frac{1}{4}) = O(n^3 \log n)$.

We can boost the coupling probability in the usual way to bound $\text{Mix}(\mathcal{M}_{\text{Gl}}, \varepsilon)$ for $\varepsilon \le 1/4$.

4.3.1. *The coupling breaks even.* Recall from Section 4.3 that $\lambda_1 = \lambda_n = \frac{1}{2}$ and $\lambda_2 = \cdots = \lambda_{n-1} = 1$.

LEMMA 4. *Suppose $(\sigma, \tau) \in S$ differs at vertex $i$. Obtain $(\sigma', \tau')$ by one step of the identity coupling. Then $\mathbf{E}[d_2(\sigma', \tau')] \le d_2(\sigma, \tau)$.*

PROOF. Recall that $n > 3$. There are three cases.

Suppose $i \in \{1, n\}$. Then $\mathbf{E}[d_2(\sigma', \tau')] - \frac{1}{2}$ is equal to

$$-\frac{2}{3n}\lambda_1 + \frac{1}{3n}\lambda_2 = 0.$$

The first term in the sum comes from the two colors which could be chosen at vertex $i$, causing coupling. The second term comes from the one bad color which could be chosen at $i$'s neighbor, causing one of the height functions to change by 2.

Suppose $i \in \{2, n-1\}$. Then $\mathbf{E}[d_2(\sigma', \tau')] - 1$ is equal to

$$\frac{2}{3n}\lambda_1 - \frac{2}{3n}\lambda_2 + \frac{1}{3n}\lambda_3 = 0.$$

Suppose $i \in \{3, \ldots, n-2\}$. Then $\mathbf{E}[d_2(\sigma', \tau')] - 1$ is equal to

$$\frac{1}{3n}\lambda_{i-1} - \frac{2}{3n}\lambda_i + \frac{1}{3n}\lambda_{i+1} = 0.$$

□

We can conclude from Lemma 4 by path-coupling that, if we take any $(\sigma, \tau) \in \Omega \times \Omega$ and produce $(\sigma', \tau')$ by one step of the identity coupling, then $\mathbf{E}[d_2(\sigma', \tau')] \le d_2(\sigma, \tau)$.

4.3.2. *Lower bounding $V$.* Let $w = \min_i \lambda_i = \frac{1}{2}$. Start with $\sigma$ and $\tau$ such that $\sigma \ne \tau$. We will identify a vertex $z$ and a color $C$ such that, if we obtain $\sigma'$ from $\sigma$ by trying $C$ at $z$ and we obtain $\tau'$ from $\tau$ by trying $C$ at $z$, then $d_2(\sigma',\tau') \le d_2(\sigma,\tau) - w$. Since $(z,C)$ is chosen with probability $1/(3n)$, we get $V = w^2/(3n)$.

Our method is this. Given $\sigma$ and $\tau$, choose $h \in \mathcal{H}(\sigma)$ and $h^* \in \mathcal{H}(\tau)$ such that $d_2(\sigma,\tau) = d(h,h^*)$. Construct $h'$ from $h$ by applying the choice $(z,C)$ (to be specified presently) in $h$ and construct $h'^*$ from $h^*$ by applying the same choice $(z,C)$ in $h^*$. We will show that $d(h',h'^*) \le d(h,h^*) - w$ so $d_2(\sigma',\tau') \le d(h',h'^*) \le d_2(\sigma,\tau) - w$.

Without loss of generality, assume that there is a vertex $v$ such that $h_v > h^*_v$. Let $m = \max_v h_v - h^*_v > 0$ and let $R = \{v | h_v - h^*_v = m\}$. By construction, $R$ is nonempty.

*Case* 1. *$R$ is the whole line.* Let $z$ be any local maximum in $h$ and let $C$ be the color that is not used at $z$ or at its neighbors in $h$. $z$ is also a local maximum in $h^*$ (since $R$ is the whole line), but $C$ is used either at $z$ or at its neighbors in $h^*$. (The unique color $C'$ that is not used either at $z$ or its neighbors in $h^*$ must be different from $C$, since $\sigma \ne \tau$.) Choose $(z,C)$. Then $h'_z = h_z - 2$. But $h'^*_z = h^*_z$. So $d(h,h^*) - d(h',h'^*) = \lambda_z$.

*Case* 2. *There is a vertex $z \in R$, all of whose neighbors are in $\overline{R}$.* Note that all edges from $z$ to $\overline{R}$ in $h$ go down (i.e., height decreases along these edges). Also, all edges from $z$ to $\overline{R}$ in $h^*$ go up. Thus, $z$ is a local maximum in $h$ and a local minimum in $h^*$. Let $C$ be the color that is not used at $z$ or at its neighbors in $h$. Choose $(z,C)$. Then $h'_z = h_z - 2$. Since $z$ is a local minimum in $h^*$, $h'^*_z \ge h^*_z$. Also, $h'_z \ge h'^*_z$ since we choose the same color in both copies. Thus, $d(h,h^*) - d(h',h'^*) \ge \lambda_z$.

*Case* 3. *There is a vertex $z \in R$ which has a neighbor $w \in \overline{R}$ and a neighbor $r \in R$.* Note that the edge from $z$ to $r$ goes the same direction (up or down) in $h$ as in $h^*$. Suppose first that it goes down. Then $z$ is a local maximum in $h$. Let $C$ be the color that is not used at $z$ or at its neighbors in $h$. Choose $(z,C)$. Then $h'_z = h_z - 2$. Also, $h'^*_z = h^*_z$ (since $z$ has a neighbor below and a neighbor above, and won't be recolored in $h^*$). Thus, $d(h,h^*) - d(h',h'^*) \ge \lambda_z$.

Suppose instead that the edge from $z$ to $r$ goes up. Then $z$ is a local minimum in $h^*$. Let $C$ be the color that is not used at $z$ or at its neighbors in $h^*$. Choose $(z,C)$. Then $h'^*_z = h^*_z + 2$ and $h'_z = h_z$ so $d(h,h^*) - d(h',h'^*) \ge \lambda_z$.

**5. Systematic scan for $q = 3$ mixes in $\Theta(n^2 \log n)$ sweeps.** As in Section 4 we consider the path $G$ with vertices 1 through $n$ with $q = 3$ colors. We consider the dynamics $\mathcal{M}_{\to}$.

5.1. *Analysis of a related Markov chain.* As in Section 4.1 the Markov chain $\mathcal{M}_{\to}$ can be associated with a Markov chain $\mathcal{M}^{\pm}_{\to}$ on $\Upsilon$. Each move of $\mathcal{M}^{\pm}_{\to}$ starts with a configuration $X \in \Upsilon$ and makes $n$ moves of the chain $\mathcal{M}_{\mathrm{Gl}}^{\pm}$ from Section 4.1 corresponding to the choices $r=1, r=2, \ldots, r=n$ (in order).

Consider the transition from configuration $X$ to configuration $X'$ corresponding to one step of $\mathcal{M}_{\rightarrow}^{\pm}$. Let $\widetilde{X}_i$ denote the label ($\pm 1$) of vertex $i$ in the intermediate configuration which is obtained after the choices $r=1, r=2, \ldots, r=i$. Then

$$\begin{aligned}
\mathbf{E}[\widetilde{X}_1] &= \tfrac{1}{3} X_1, \\
\mathbf{E}[\widetilde{X}_i] &= \tfrac{2}{3} X_i + \tfrac{1}{3} \mathbf{E}[\widetilde{X}_{i-1}], \qquad && i = 2, \ldots, n-1, \\
\mathbf{E}[X'_i] &= \tfrac{2}{3} \mathbf{E}[\widetilde{X}_i] + \tfrac{1}{3} X_{i+1}, \qquad && i = 1, \ldots, n-2, \\
\mathbf{E}[X'_{n-1}] &= \tfrac{1}{3} \mathbf{E}[\widetilde{X}_{n-1}].
\end{aligned}$$

Solving these gives

$$\begin{aligned}
\mathbf{E}[X'_1] &= \tfrac{2}{9} X_1 + \tfrac{1}{3} X_2, \\
\mathbf{E}[X'_i] &= \frac{2}{3^{i+1}} X_1 + \sum_{j=2}^{i} \frac{4}{3^{i+2-j}} X_j + \frac{1}{3} X_{i+1}, \qquad i = 2, \ldots, n-2, \\
\mathbf{E}[X'_{n-1}] &= \frac{1}{3^n} X_1 + \sum_{j=2}^{n-2} \frac{2}{3^{n+1-j}} X_j + \frac{2}{9} X_{n-1}.
\end{aligned}$$

So the matrix

$$A = \begin{bmatrix}
\frac{2}{9} & \frac{1}{3} & 0 & & & & \\
\frac{2}{27} & \frac{4}{9} & \frac{1}{3} & 0 & & & \\
\frac{2}{81} & \frac{4}{27} & \frac{4}{9} & \frac{1}{3} & 0 & & \\
\vdots & \vdots & \vdots & \ddots & \ddots & \ddots & \\
\vdots & \vdots & \vdots & & \ddots & \ddots & 0 \\
\frac{2}{3^{n-1}} & \frac{4}{3^{n-2}} & \frac{4}{3^{n-3}} & \frac{4}{3^{n-4}} & \cdots & \frac{4}{9} & \frac{1}{3} \\
\frac{1}{3^n} & \frac{2}{3^{n-1}} & \frac{2}{3^{n-2}} & \frac{2}{3^{n-3}} & \cdots & \frac{2}{27} & \frac{2}{9}
\end{bmatrix}.$$

Here $A$ is not symmetric, but is nonnegative and irreducible, so has a positive eigenvector $w$ corresponding to its (real) largest eigenvalue $\lambda$. Now $w, \lambda$ satisfy the equations

$$\lambda w_1 = \frac{2}{9} w_1 + \sum_{j=2}^{n-2} \frac{2}{3^{j+1}} w_j + \frac{1}{3^n} w_{n-1},$$

$$\lambda w_i = \frac{1}{3} w_{i-1} + \sum_{j=i}^{n-2} \frac{4}{3^{j-i+2}} w_j + \frac{2}{3^{n-i+1}} w_{n-1}, \qquad i = 2, \ldots, n-2,$$

$$\lambda w_{n-1} = \tfrac{1}{3} w_{n-2} + \tfrac{2}{9} w_{n-1}.$$

These can be simplified by subtracting one-third of the $(i+1)$st equation from the $i$th for $i = 2, \ldots, n-2$, and one-sixth the second from the first, giving

$$\lambda w_2 - (6\lambda - 1) w_1 = 0, \tag{9}$$

$$\lambda w_{i+1} - (3\lambda - 1) w_i + w_{i-1} = 0, \qquad i = 2, \ldots, n-2, \tag{10}$$

$$-3 w_{n-2} + (9\lambda - 2) w_{n-1} = 0. \tag{11}$$

If $\lambda$ is close to 1, the form of (10) suggests a slightly damped harmonic oscillation, so we will try a solution of the form $w_i = c_n e^{\gamma i} \sin(\alpha i + \beta)$, where $c_n > 0$ is a constant, depending on $n$, that can be chosen later. Substituting this in (10) and equating coefficients of $\sin(\alpha i + \beta)$, $\cos(\alpha i + \beta)$ gives

$$\lambda = e^{-2\gamma} \quad \text{and} \quad \cos\alpha = (3e^{-\gamma} - e^{\gamma})/2, \qquad \text{that is, } e^{\gamma} = \sqrt{3 + \cos^2\alpha} - \cos\alpha. \tag{12}$$

[The second of these follows from $\sin(x+y) = \sin x \cos y + \cos x \sin y$ and $\sin(x-y) = \sin x \cos y - \cos x \sin y$ and the third used the quadratic formula with the choice $\cos\alpha \geq 0$.] Then (9) and (11) give

$$\frac{\sin(2\alpha + \beta)}{\sin(\alpha + \beta)} = \cos\alpha + \cot(\alpha + \beta) \sin\alpha = 6e^{-\gamma} - e^{\gamma} \tag{13}$$

and

$$\frac{\sin((n-2)\alpha + \beta)}{\sin((n-1)\alpha + \beta)} = \cos\alpha - \cot((n-1)\alpha + \beta) \sin\alpha = \frac{9e^{-\gamma} - 2e^{\gamma}}{3}. \tag{14}$$

Using (12) to eliminate $\gamma$ in (13) and (14) gives

$$\tan(\alpha + \beta) = \frac{\sin\alpha}{2\cos\alpha + \sqrt{3 + \cos^2\alpha}}$$

and

$$\tan((n-1)\alpha + \beta) = \frac{-3\sin\alpha}{2\cos\alpha + \sqrt{3 + \cos^2\alpha}},$$

implying

$$\tan((n-1)\alpha + \beta) = -3\tan(\alpha + \beta) \quad \text{and} \quad \tan\alpha = 4\tan(\alpha + \beta)/(1 - 3\tan^2(\alpha + \beta)). \tag{15}$$

To see the second of these equalities, note that the left-hand equation on the previous line is equivalent to

$$\tan(\alpha+\beta) = \frac{\sin\alpha}{2\cos\alpha + \sqrt{4\cos^2\alpha + 3\sin^2\alpha}}$$

using $\cos^2\alpha + \sin^2\alpha = 1$. But this is equal to

$$\frac{\tan\alpha}{2+\sqrt{4+3\tan^2\alpha}}.$$

Now solve this for $\tan\alpha$. The equalities in (15) imply

$$\pi-(n-1)\alpha-\beta = \arctan(3\tan(\alpha+\beta)),$$
$$\alpha = \alpha+\beta+\arctan(3\tan(\alpha+\beta)).$$

The first of these uses $\tan(\pi-x) = -\tan(x)$ and the second uses $\tan(x+y) = (\tan x + \tan y)/(1-\tan x \tan y)$, with $x = \alpha+\beta$ and $y = \arctan(3\tan(\alpha+\beta))$. So finally we have

$$\alpha = \frac{\pi}{n-1},$$
$$\tan\beta = -3\tan\left(\beta+\frac{\pi}{n-1}\right),$$
$$e^\gamma = \sqrt{3+\cos^2\left(\frac{\pi}{n-1}\right)} - \cos\left(\frac{\pi}{n-1}\right).$$

Note that $\beta$ is the solution of a trigonometric equation, but it is easily checked that $-\pi/(n-1) < \beta < 0$. Hence, $w > 0$, corresponding to the largest eigenvalue $\lambda$ of $A$. Asymptotically, we have

$$\alpha \sim \pi/n, \qquad \beta \sim -3\pi/4n, \qquad \gamma \sim \pi^2/4n^2, \qquad \text{so } 1-\lambda \sim \pi^2/2n^2. \tag{16}$$

We also need to set $c_n \sim 4n/\pi$ to achieve $w_i \geq 1$, for all $0 < i < n$. If we take $X(0)$ to be the all 1's vector, then it is easy to check that $\Phi_0 \sim 8n^2/\pi^2$.

Next we need to estimate $\rho$, the bound on the variance of $\Phi_t$, given $\Phi_{t-1}$. In the case of Glauber dynamics, the range of $\Phi_t$ was $O(1)$, which provided a crude bound $\rho = O(1)$. For scan, however, the range of possible values of $\Phi_t$ is $O(n)$, which yields only $\rho = O(n^2)$: too weak for our purposes. Intuitively, however, since $\Phi_t$ is, roughly speaking, a sum of $n$ nearly independent r.v.'s each of variance $O(1)$, the variance of $\Phi_t$ ought to be $O(n)$. This is indeed the case. In fact, we prove something stronger in the form of a large deviation result for $\Phi_t$. Before doing that, let's complete the remainder of the proof. Assuming $\rho = O(n)$, we have $\nu = O(n^3)$. Now, from (3), the lower bound is $\mathrm{Mix}(\mathcal{M}^{\pm}_{\rightarrow}, \frac{1}{2}) \geq \pi^{-2}n^2\ln n - O(n^2)$. Since the sweep is also monotone, (5) gives the upper bound $\mathrm{Mix}(\mathcal{M}^{\pm}_{\rightarrow}, \varepsilon) \leq 4\pi^{-2}n^2\ln n + \frac{2n^2}{\pi^2}\ln\varepsilon^{-1} + O(n^2)$. It may be observed that the bounds for $\mathcal{M}_{\mathrm{Gl}}^{\pm}$ are both about $n$ times these quantities, so there is no evidence that the scan gives a significant speed-up. However, there will be a considerable saving in random number generation.

It only remains to show $\rho = O(n)$. Recall that $\rho$ is an upper bound on $\mathbf{E}[\mathrm{var}(\Phi_t|\Phi_{t-1})]$, where $\Phi_t = wX(t)$. Also, for $i \in \{1, \ldots, n-1\}$, $w_i = (4n/\pi)\exp(\gamma i)\sin(\alpha i + \beta)$, where $\gamma$, $\alpha$ and $\beta$ are given asymptotically in (16). Let $w_0 = w_n = 0$. Our first observation, which is similar to the one used in the analysis of Glauber dynamics, is

$$\max_{1 \le i \le n} |w_i - w_{i-1}| = O(1). \tag{17}$$

To see that (17) holds, first note that $1 \le \exp(\gamma i) = 1 + O(1/n)$. Using the series expansion of sine, we find that $w_1 = O(1)$ and $w_{n-1} = O(1)$. Now, for $i \in \{2, \ldots, n-1\}$, note that

$$\begin{aligned} w_i - w_{i-1} &\le \frac{4n}{\pi}(1 + O(1/n))\sin(\alpha i + \beta) - \frac{4n}{\pi}\sin(\alpha(i-1) + \beta) \\ &\le O(1)\sin(\alpha i + \beta) + \frac{4n}{\pi}(\sin(\alpha i + \beta) - \sin(\alpha(i-1) + \beta)). \end{aligned}$$

The first term is $O(1)$ because sine is bounded. Since the derivative of sine is at most 1 (in absolute value), the difference between the two sines in the second term is at most $\alpha$, so the second term is also $O(1)$.

Let $\omega_1, \omega_2, \ldots, \omega_n$ denote the sequence of swap/no-swap decisions made by systematic scan in transforming $X(t-1)$ to $X(t)$. That is, $\omega_1$ is the indicator r.v. for the event that the sign of position 1 is flipped, $\omega_i$ (for $i \in \{2, \ldots, n-1\}$) is the indicator r.v. for the event that positions $i-1$ and $i$ are exchanged, and $\omega_n$ is the indicator r.v. for the event that the sign of position $n$ is flipped. The $\omega_i$'s are independent Bernoulli random variables with parameter $1/3$. Given $X(t-1)$, the configuration $X(t)$ is a r.v. in $\omega_1, \omega_2, \ldots, \omega_n$. Let $\omega_{n+1} = 0$. Consider the Doob martingale $Z_0, Z_1, \ldots, Z_n$ obtained by revealing the swap/no-swap decisions in sequence:

$$\begin{aligned} &Z_0 = \mathbf{E}[wX(t)], \qquad Z_1 = \mathbf{E}[wX(t)|\omega_1], \\ &Z_2 = \mathbf{E}[wX(t)|\omega_1, \omega_2], \ldots, Z_n = \mathbf{E}[wX(t)|\omega_1, \omega_2, \ldots, \omega_n]. \end{aligned}$$

All of $Z_0, \ldots, Z_n$ are conditioned on $X(t-1)$. Notice that $Z_0 = \mathbf{E}[\Phi_t|X(t-1)]$ and $Z_n = \Phi_t$. We will show below that $|Z_{i-1} - Z_i| = O(1)$, for all $1 \le i \le n$. It follows from the Azuma–Hoeffding inequality [4] (see also [21], Chapter 2.4) that

$$\Pr(|\Phi_t - \mathbf{E}[\Phi_t]| > h\sqrt{n}) = \Pr(|Z_n - Z_0| > h\sqrt{n}) \le \exp(-\Omega(h^2)).$$

Let $C = \max_i |Z_{i-1} - Z_i|$. Then we get $\rho = O(n)$ since

$$\mathbf{E}[\mathrm{var}(\Phi_t|\Phi_{t-1})]$$

$$
\begin{aligned}
&= \sum_{\xi} \Pr(X(t-1)=\xi)\mathbf{E}[(\Phi_t - \mathbf{E}[\Phi_t|X(t-1)=\xi])^2|X(t-1)=\xi] \\
&\le \max_{\xi} \mathbf{E}[(\Phi_t - \mathbf{E}[\Phi_t|X(t-1)=\xi])^2|X(t-1)=\xi] \\
&= \max_{\xi} \mathbf{E}[(Z_n - Z_0)^2|X(t-1)=\xi] \\
&\le n \max_{\xi} \sum_{h=1}^{\lceil C\sqrt{n}\rceil} h^2 \Pr\left(|Z_n - Z_0| \in ((h-1)\sqrt{n}, h\sqrt{n})|X(t-1)=\xi\right) \\
&\le n \sum_{h=0}^{\lceil C\sqrt{n}\rceil - 1} (h+1)^2 \exp(-\Omega(h^2)) = O(n).
\end{aligned}
$$

Finally, we must argue that $|Z_{i-1} - Z_i| = O(1)$. First, note that

$$
\begin{aligned}
|Z_1 - Z_0| &= |\mathbf{E}[wX(t)|\omega_1] - \mathbf{E}[wX(t)]| \\
&\le |\mathbf{E}[wX(t)|\omega_1 = 1] - \mathbf{E}[wX(t)|\omega_1 = 0]|,
\end{aligned}
$$

and the right-hand side is at most

$$
\Bigg| \sum_{k=0}^{n-1} \Pr\left((\omega_2, \ldots, \omega_{k+2}) = (1, \ldots, 1, 0)\right) \times \left(\mathbf{E}[wX(t)|1, 1, \ldots, 1, 0] - \mathbf{E}[wX(t)|0, 1, \ldots, 1, 0]\right) \Bigg|,
$$

where the conditioning specifies the values of $\omega_1, \ldots, \omega_{k+2}$. The relevant probability is at most $3^{-k}$. To get an upper bound, we move the absolute value inside the summation and maximise over $\omega_{k+3}, \ldots, \omega_n$, obtaining

$$
\begin{aligned}
|Z_1 - Z_0| \le \sum_{k=0}^{n-1} 3^{-k} \max_{\omega_{k+3}, \ldots, \omega_n} &|\mathbf{E}[wX(t)|1, 1, \ldots, 1, 0, \omega_{k+3}, \ldots, \omega_n] \\
&- \mathbf{E}[wX(t)|0, 1, \ldots, 1, 0, \omega_{k+3}, \ldots, \omega_n]| \\
= \sum_{k=0}^{n-1} 3^{-k} 2 w_{k+1},
\end{aligned}
$$

since the difference in sign propagates to position $k+1$ and then stops. By (17), this is $O(1)$. Similarly, $|Z_n - Z_{n-1}| = O(1)$. Now consider $i \in \{2, \ldots, n-1\}$. Mimicking the analysis above, we find that $|Z_i - Z_{i-1}|$ is at most

$$
\sum_{k=0}^{n-i} 3^{-k} \max_{\omega_{i+k+2}, \ldots, \omega_n} T,
$$

where $T$ is the absolute value of

$$\mathbf{E}[wX(t)|\omega_1,\dots,\omega_{i-1},1,1,\dots,1,0,\omega_{i+k+2},\dots,\omega_n] \\ -\mathbf{E}[wX(t)|\omega_1,\dots,\omega_{i-1},0,1,\dots,1,0,\omega_{i+k+2},\dots,\omega_n],$$

which is at most $2|w_{i+k}-w_{i-1}|$, so

$$|Z_i - Z_{i-1}| \le \sum_{k=0}^{n-i} 3^{-k} 2|w_{i+k}-w_{i-1}|,$$

which is $O(1)$ by (17).

5.2. *A lower bound for systematic scan.* We will use distance measures $d_1$ and $d_2$ from Section 4.2. Using distance measure $d_1$, the lower bound from Section 5.1 applies directly to $\mathcal{M}_\rightarrow$. Thus, we have the following:

THEOREM 5. *Let $G$ be the $n$-vertex path, and let $q=3$. Then a lower bound on the mixing time of the Markov chain $\mathcal{M}_\rightarrow$ on the state space $\Omega$ is given by* $\mathrm{Mix}(\mathcal{M}_\rightarrow,\frac{1}{2}) \ge \pi^{-2}n^2\ln n - O(n)$.

5.3. *An upper bound for systematic scan.* As in Section 4.3, we find that, for some $t' = O(n^2\log n)$, we will have $d_1(\sigma(t'),\tau(t'))=0$ with probability at least $\frac{39}{40}$.

Now choose the following weights. Let $\lambda_1 = \frac{1}{4}$, $\lambda_2 = \cdots = \lambda_{n-1} = 1$ and $\lambda_n = \frac{3}{4}$.

We now use path coupling on pairs $(\sigma(0),\tau(0)) \in S$. Start with such a pair, run $t'$ steps to get $(\sigma(t'),\tau(t'))$. With probability at least 39/40, $d_1(\sigma(t'),\tau(t'))=0$. Either $\sigma(t')=\tau(t')$ or $d_2(\sigma(t'),\tau(t'))=n-1$. Suppose the latter. We now carry on from $(\sigma(t'),\tau(t'))$ using the identity coupling. We show in Section 5.3.1 that, if we take any $(\sigma,\tau)\in\Omega\times\Omega$ and produce $(\sigma',\tau')$ by one scan using the identity coupling, then $\mathbf{E}[d_2(\sigma',\tau')] \le d_2(\sigma,\tau)$. Thus, $D_t = d_2(\sigma(t'+t),\tau(t'+t))$ is a super-martingale with $D_0 = n-1$. In Section 5.3.2, we define $V = 1/27$ and show that, for all $t$ and all values of $D_t$ other than 0, $\mathbf{E}[(D_{t+1}-D_t)^2|D_t] \ge V$. Let $B = 10n$, and let $T$ be the first time at which either (a) $D_t = 0$ (i.e., coupling occurs), or (b) $D_t \ge B$. Note that $T$ is a stopping time. Define $Z_t = (B-D_t)^2 - Vt$, and observe, as in Section 4.3, that $Z_{t\wedge T}$ is a sub-martingale. Let $p$ be the probability that (a) occurs. As in Section 4.3, applying the optional stopping theorem to $D_T$ gives $p \ge \frac{9}{10}$. Also, as before,

$$pB^2 + (1-p)E[(B-D_T)^2|D_T \ge B] - VE[T] > 0.$$

Since $|D_t - D_{t-1}| \le 2n$, $(1-p)E[(B-D_T)^2|D_T\ge B] \le 4n^2 < pB^2$ so $\mathbf{E}[T] \le (2pB^2)/V$. Conditioning on (a) occurring, it follows that $\mathbf{E}[T|D_T=0] \le$

$2B^2/V$. Hence, $\Pr(T > 20B^2/V \mid D_T = 0) \le \frac{1}{10}$. So, if we now run the identity coupling for $20B^2/V = O(n^2)$ steps, then $\sigma$ and $\tau$ will fail to couple with probability at most $\frac{1}{40} + 2 \times \frac{1}{10} < \frac{1}{4}$. Thus, we have shown the following:

THEOREM 6. *Let $G$ be the $n$-vertex path, and let $q = 3$. Consider the Markov chain $\mathcal{M}_{\to}$ on the state space $\Omega$. Then $\mathrm{Mix}(\mathcal{M}_{\to}, \frac{1}{4}) = O(n^2 \log n)$.*

We can bound $\mathrm{Mix}(\mathcal{M}_{\to}, \varepsilon)$ for $\varepsilon \le 1/4$ by boosting the coupling probability in the usual way.

5.3.1. *The coupling breaks even.* Recall that the vertices of the path $G$ are labeled $1, \dots, n$ going from left to right.

LEMMA 7. *Suppose that $\sigma$ and $\tau$ differ at vertex $i < n$ and agree to the right of vertex $i$. Obtain $\sigma'$ and $\tau'$ by scanning left to right, starting at vertex $i+1$, doing the identity coupling. Then*

$$\mathbf{E}[d_2(\sigma', \tau')] - d_2(\sigma, \tau) \le \tfrac{1}{2}.$$

PROOF. Choose $h \in \mathcal{H}(\sigma)$ and $h^* \in \mathcal{H}(\tau)$ such that $d_2(\sigma, \tau) = d(h, h^*)$. Let $h', h'^*$ be the transformed height functions produced by the scan. If $v = i + \ell$ for $\ell \in \{1, \dots, n-i-1\}$, then the probability that vertex $v$ is changed by the coupling is $(\frac{1}{3})^\ell$. If there is a change, then one of the height functions changes by 2, so the change in $d(h', h'^*)$ is 1. For $v = n$, the probability that $v$ changes is $(\frac{1}{3})^{n-i-1} \frac{2}{3}$. The change to $d(h', h'^*)$ in this case is $\frac{3}{4}$. Thus,

$$\begin{aligned}
\mathbf{E}[d_2(\sigma', \tau')] &\le \mathbf{E}[d(h', h'^*)] \\
&= d(h, h^*) + \sum_{\ell=1}^{n-i-1} \left(\frac{1}{3}\right)^\ell + \frac{2}{3}\frac{3}{4}\left(\frac{1}{3}\right)^{n-i-1} \\
&= d(h, h^*) + \frac{1/3 - (1/3)^{n-i}}{1 - 1/3} + \frac{1}{2}\left(\frac{1}{3}\right)^{n-i-1} \\
&= d(h, h^*) + \frac{1}{2}.
\end{aligned}$$

□

LEMMA 8. *Suppose that $\sigma$ and $\tau$ differ only at vertex 1. Obtain $\sigma'$ and $\tau'$ by scanning left to right, starting at vertex 1, doing the identity coupling. Then*

$$\mathbf{E}[d_2(\sigma', \tau')] \le \tfrac{1}{4}.$$

PROOF. With probability $\frac{2}{3}$, the first vertex agrees, so $\sigma' = \tau'$. With probability $\frac{1}{3}$, the first vertex is left unchanged. Thus (using Lemma 7), $\mathbf{E}[d_2(\sigma', \tau')] \le \frac{1}{3}(\frac{1}{4} + \frac{1}{2})$. □

LEMMA 9. *Suppose that $\sigma$ and $\tau$ differ only at vertex 2. Obtain $\sigma'$ and $\tau'$ by scanning left to right, starting at vertex 1, doing the identity coupling. Then*

$$\mathbf{E}[d_2(\sigma',\tau')] \le 1.$$

PROOF. Say that $\sigma$ starts 2 0 2 and $\tau$ starts 2 1 2. Consider the coupling of first vertex:

- With probability $\frac{2}{3}$:
  The first vertex is made to disagree, for example, $\sigma$ now starts 2 0 2 but $\tau$ starts 0 1 2.
  The coupling of the second vertex in this case is as follows:
  - With probability $\frac{1}{3}$:
    The second vertex is made to agree, for example, both become 1. In this case, $\mathbf{E}[d_2(\sigma',\tau')] = \frac{1}{4}$.
  - With probability $\frac{2}{3}$:
    The second vertex is unchanged. In this case, $\mathbf{E}[d_2(\sigma',\tau')] \le \frac{1}{4} + 1 + \frac{1}{2} = \frac{7}{4}$.
- With probability $\frac{1}{3}$:
  The first vertex is unchanged. By analogy to the proof of Lemma 8, $\mathbf{E}[d_2(\sigma',\tau')] \le \frac{1}{3}(1+\frac{1}{2}) = \frac{1}{2}$.

Adding it all up, $\mathbf{E}[d_2(\sigma',\tau')] \le \frac{2}{3}(\frac{1}{3}\cdot\frac{1}{4} + \frac{2}{3}\cdot\frac{7}{4}) + \frac{1}{3}\cdot\frac{1}{2} = 1$. □

LEMMA 10. *Let $2 < i < n$. Suppose that $\sigma$ and $\tau$ differ only at vertex $i$. Obtain $\sigma'$ and $\tau'$ by scanning left to right, starting at vertex 1, doing the identity coupling. Then*

$$\mathbf{E}[d_2(\sigma',\tau')] \le 1.$$

PROOF. Say that vertices $i-2,\dots,i+1$ of $\sigma$ are 1 2 0 2 and of $\tau$ are 1 2 1 2.
Consider the coupling of vertex $i-1$:

- With probability $\frac{1}{3}$:
  Vertex $i-1$ is made to disagree, so $\sigma$ now starts 1 2 0 2 but $\tau$ starts 1 0 1 2. By analogy with the proof of Lemma 9, $\mathbf{E}[d_2(\sigma',\tau')] \le \frac{1}{3}\cdot 1 + \frac{2}{3}(1+1+\frac{1}{2}) = 2$.
- With probability $\frac{2}{3}$:
  Vertex $i-1$ is unchanged. By analogy with the proof of Lemma 8, $\mathbf{E}[d_2(\sigma',\tau')] \le \frac{1}{3}(1+\frac{1}{2}) = \frac{1}{2}$.

Adding it all up, $\mathbf{E}[d_2(\sigma',\tau')] \le \frac{1}{3}\cdot 2 + \frac{2}{3}\cdot\frac{1}{2} = 1$. □

LEMMA 11. *Suppose that $\sigma$ and $\tau$ differ only at vertex $n$. Obtain $\sigma'$ and $\tau'$ by scanning left to right, starting at vertex 1, doing the identity coupling. Then*

$$\mathbf{E}[d_2(\sigma', \tau')] \le \tfrac{3}{4}.$$

PROOF. Say that vertices $n-2, \ldots, n$ of $\sigma$ are $0\,2\,0$ and of $\tau$ are $0\,2\,1$. Consider the coupling of vertex $n-1$:

- With probability $\frac{1}{3}$:
  Vertex $n-1$ is made to disagree, so $\sigma$ now ends with $0\,1\,0$ and $\tau$ with $0\,2\,1$.
  The coupling of vertex $n$ (the last) in this case is as follows:
  - With probability $\frac{1}{3}$:
    The last vertex is made 0, with resulting cost 1.
  - With probability $\frac{1}{3}$:
    The last vertex is unchanged, with resulting cost $1 + \frac{3}{4} = \frac{7}{4}$.
  - With probability $\frac{1}{3}$:
    The last vertex becomes 2 in $\sigma$, with resulting cost $1 + 2 \cdot \frac{3}{4} = \frac{5}{2}$. (The claimed final cost is witnessed by the sequence of transitions $0\,1\,2 \to 0\,1\,0 \to 0\,2\,0 \to 0\,2\,1$.)
- With probability $\frac{2}{3}$:
  Vertex $n-1$ is unchanged. Now with probability $\frac{2}{3}$, vertex $n$ will agree and with probability $\frac{1}{3}$, it will be unchanged. Thus, $\mathbf{E}[d_2(\sigma', \tau')] \le \frac{1}{3} \cdot \frac{3}{4} = \frac{1}{4}$.

Adding it all up, $\mathbf{E}[d_2(\sigma', \tau')] \le \frac{1}{3}(\frac{1}{3} \cdot 1 + \frac{1}{3} \cdot \frac{7}{4} + \frac{1}{3} \cdot \frac{5}{2}) + \frac{2}{3} \cdot \frac{1}{4} = \frac{3}{4}$. □

Lemmas 8, 9, 10 and 11 show that if $(\sigma, \tau) \in S$ and we obtain $\sigma'$ and $\tau'$ by scanning left to right, starting at vertex 1, doing the identity coupling, then $\mathbf{E}[d_2(\sigma', \tau')] \le d_2(\sigma, \tau)$. By path coupling, we find that if we take any $(\sigma, \tau) \in \Omega \times \Omega$ and we produce $(\sigma', \tau')$ by one scan using the identity coupling, then $\mathbf{E}[d_2(\sigma', \tau')] \le d_2(\sigma, \tau)$.

5.3.2. *The coupling has enough variance* (*lower bounding* $V$). Recall that $w = \min_i \lambda_i$. Suppose $\sigma \ne \tau$. In Section 4.3.2, we considered several cases. For each case, we identified a vertex $z$ and a color $C$ such that, if we obtain $\sigma'$ from $\sigma$ by trying $C$ at $z$ and we obtain $\tau'$ from $\tau$ by trying $C$ at $z$, then $d_2(\sigma', \tau') \le d_2(\sigma, \tau) - w$.

In this section we reconsider each case. Obtain $\sigma^*$ and $\tau^*$ by scanning $\sigma$ and $\tau$ left to right, using the identity coupling. For each case in Section 4.3.2, we prove the following:

- If $z > 1$, then there is a color $c_\ell$ (depending only on $\sigma_{z-1}$, $\sigma_z$, $\tau_{z-1}$ and $\tau_z$) such that choosing color $c_\ell$ for vertex $z-1$ ensures $\sigma^*_{z-1} = \sigma_{z-1}$ and $\tau^*_{z-1} = \tau_{z-1}$. Actually, it is easy to see that $c_\ell$ exists—just take any color in $\{\sigma_{z-1}, \sigma_z\} \cap \{\tau_{z-1}, \tau_z\}$.

TABLE 1

| **Case** | $\boldsymbol{\sigma_{z-1}\cdots\sigma_{z+1}}$ | $\boldsymbol{\tau_{z-1}\cdots\tau_{z+1}}$ | $\boldsymbol{C'}$ | $\boldsymbol{c_r(0)}$ | $\boldsymbol{c_r(1)}$ | $\boldsymbol{c_r(2)}$ |
|---|---|---|---|---|---|---|
| 1 | 010 | 121 | 1 | 0 | 1 | 2 |
| 1 | 010 | 202 | 0 | 0 | 1 | 2 |
| 2 | 010 | 101 | 0 | 0 | 0 | 2 |
| 2 | 010 | 212 | 1 | 0 | 1 | 2 |
| 2 | 010 | 020 | 0 | 0 | 0 | 2 |

- Suppose $z < n$. For any color $c$, there is a color $c_r(c)$ (depending only on $\sigma_{z-1}$, $\sigma_z$, $\sigma_{z+1}$, $\tau_{z-1}$, $\tau_z$, $\tau_{z+1}$ and $c$) such that if we choose $c_\ell$ for vertex $z-1$, $c$ for vertex $z$ and $c_r(c)$ for vertex $z+1$, then $\sigma^*_{z+1} = \sigma_{z+1}$ and $\tau^*_{z+1} = \tau_{z+1}$.
- There is a color $C'$ such that, if we obtain $\sigma'$ from $\sigma$ by trying $C'$ at $z$ and we obtain $\tau'$ from $\tau$ by trying $C'$ at $z$, then $\sigma'_z = \sigma_z$ and $\tau'_z = \tau_z$.

This is enough to establish $V = 1/27$. We will consider the event that $c_\ell$ is chosen for $z-1$ and, whatever color, $c$, is chosen for $z$, $c_r(c)$ is chosen for $z+1$. This event occurs with probability 1/9. Conditioned on the fact that this event occurs, we can choose the color $c$ for vertex $z$ *after* choosing all other colors. That is, the choice of $c$ is independent of the rest of the scan. Let $\sigma^\dagger$ and $\tau^\dagger$ be random variables defined by a left to right scan of $\sigma$ and $\tau$, which uses $c_\ell$ at $z-1$ and $c_r(c)$ at $z+1$ and misses out the re-coloring at $z$.

If $|d_2(\sigma^\dagger, \tau^\dagger) - d_2(\sigma, \tau)| \geq w/2$, then we choose color $C'$ for vertex $z$ so $\sigma^* = \sigma^\dagger$ and $\tau^* = \tau^\dagger$. Otherwise, we choose color $C$ for vertex $v$ so $d_2(\sigma^*, \tau^*) \leq d_2(\sigma^\dagger, \tau^\dagger) - w$. Either way, we get $|d_2(\sigma^*, \tau^*) - d_2(\sigma, \tau)| \geq w/2$.

Cases 1 and 2 from Section 4.3.2 are in Table 1.

In Case 3, say $\sigma_{z-1}\cdots\sigma_{z+1} = 010$ and $\tau_{z-1}\cdots\tau_{z+1}$ is monotonic. Then $C'$ is any color in $\{0, 1\}$, $c_r(2)$ is any color in $\{0, 2\} \cap \{\tau_z, \tau_{z+1}\}$ and, for any $i \neq 2$, $c_r(i)$ is any color in $\{0, 1\} \cap \{\tau_z, \tau_{z+1}\}$.

The case where $\tau_{z-1}\cdots\tau_{z+1} = 101$ and $\sigma_{z-1}\cdots\sigma_{z+1}$ is monotonic is similar.

## 6. Optimal mixing of Glauber and scan when $q = 4$.

6.1. *Distance measures.* In this section $G$ is the $n$-vertex path. We take the state space to be $\Omega^+$ (i.e., all colorings, whether proper or not). The results that we get by analyzing our Markov chains on state space $\Omega^+$ also apply to the same chains with state space $\Omega$—this is because the chains do not make transitions from states in $\Omega$ to states outside of $\Omega$. (Thus, the stationary distribution is uniform on $\Omega$ – states in $\Omega^+ \setminus \Omega$ have zero measure.) We ought to note that, on the extra states in $\Omega^+ \setminus \Omega$, what we are calling a "Metropolis" update does not strictly fit the official definition. For example, with a natural definition of the "energy" of a coloring, and using the usual Metropolis filter, the transition $\dots 0\,0\,1 \dots \to \dots 0\,1\,1 \dots$ would occur with positive probability. Nevertheless, we disallow this transition because of the adjacent color 1 vertices in the final state. However, our version of Metropolis agrees with the usual one on the significant part of the state space, namely, $\Omega$.

6.2. *Glauber with $q = 4$: $O(n \log n)$ updates suffice.* We'll use Theorem 2.2 of [13].

Suppose $(\sigma, \tau) \in S$ differ on vertex $i$. Construct $(\sigma', \tau')$ from $(\sigma, \tau)$ by using the following coupling. Choose the same vertex $v$ to recolor in $\sigma$ and in $\tau$. Choose the same color in both copies unless $v \in \{i-1, i+1\}$. In that case, choose color $\sigma_i$ in one copy, while choosing $\tau_i$ in the other (and choose the same color otherwise).

We first show that the value $\beta$ in Theorem 2.2 of [13] is 1. That is, we show that $\mathbf{E}[\mathrm{Ham}(\sigma', \tau')] \le 1$. Consider the choices made in $\sigma$. If we choose vertex $i-1$ and color $\tau_i$, then the Hamming distance might go up by 1. Similarly, if we choose vertex $i+1$ and color $\tau_i$, then the Hamming distance might go up by 1. If we choose vertex $i$ and any of the (at least two) colors not in $\{\sigma_{i-1}, \sigma_{i+1}\}$, then the Hamming distance goes down by 1. These are the only choices which can cause the distance to change.

Now consider a multi-step coupling from $(\sigma, \tau)$. Assume for now that $i \in \{3, \dots, n-2\}$, so there are at least 2 vertices to the left of vertex $i$ and at least 2 vertices to the right of vertex $i$. The other cases are easier and we will consider them later. Let $c$ be a color which is not in $\{\sigma_i, \tau_i\}$ (there are two such colors, but choose an arbitrary one and call it $c$). Let $\Psi$ be the set containing the following 6 choices (in $\sigma$): choose $i$ with any color, choose $i-1$ with $\tau_i$, or choose $i+1$ with $\tau_i$. Let the stopping time $T$ be the first time a choice from $\Psi$ is made. [I.e., a choice from $\Psi$ is made in the transition from $(\sigma(T-1), \tau(T-1))$ to $(\sigma(T), \tau(T))$, where $(\sigma(0), \tau(0)) = (\sigma, \tau)$.]

Let $\Xi$ be the set containing the following 14 choices (in $\sigma$): Choose $i-1$ with any color besides $\tau_i$. Choose $i+1$ with any color besides $\tau_i$. Choose $i-2$ with any color. Choose $i+2$ with any color. Let $C$ be the set containing all $4n-20$ choices that are not in $\Psi$ or $\Xi$. Let $z_1, \dots, z_t$ denote the choices made (in $\sigma$) in the transitions $(\sigma(0), \tau(0)), (\sigma(1), \tau(1)), \dots, (\sigma(t), \tau(t))$. We will say that the sequence $z_1, \dots, z_t$ is *good* if the only choices in $\Xi \cup \Psi$ are the following:

- for some $t_1 \in [1, t]$, $z_{t_1}$ consists of vertex $i-2$ along with the "smallest" color (e.g., smallest numerically) that is not in $\{c, \sigma_{i-3}(t_1 - 1), \sigma_{i-1}\}$, and

- for some $t_2 \in [t_1+1, t]$, $z_{t_2}$ consists of vertex $i+2$ along with the smallest color that is not in $\{c, \sigma_{i+3}(t_1-1), \sigma_{i+1}\}$, and
- for some $t_3 \in [t_2+1, t]$, $z_{t_3} = (i-1, c)$,
- for some $t_4 \in [t_3+1, t]$, $z_{t_4} = (i+1, c)$.

Denote by $\mathcal{G}$ the event that $z_1, \ldots, z_{t-1}$ is good. Now,

$$\Pr(\mathcal{G} | T = t) = \binom{t-1}{4} \left(\frac{1}{4n-6}\right)^4 \left(\frac{4n-20}{4n-6}\right)^{t-5}. \tag{18}$$

Let $\alpha$ be any positive constant which is at most 1/6. Let $\delta$ be a positive constant, independent of $n$, such that, for all $t \in [\alpha n, n]$, the expression in (18) is at least $\delta$. Now

$$\mathbf{E}[\,\mathrm{Ham}(\sigma(T), \tau(T)) | T = t \text{ and } \mathcal{G}] \le \frac{3\times 0 + 1\times 1 + 2\times 2}{6} = \frac{5}{6}.$$

In particular, $\mathrm{Ham}(\sigma(T), \tau(T)) = 1$ if $z_t = (i, c)$ and $\mathrm{Ham}(\sigma(T), \tau(T)) = 0$ if $z_t$ consists of vertex $i$ with some other color. Otherwise, $\mathrm{Ham}(\sigma(T), \tau(T)) \le 2$. Similarly,

$$\mathbf{E}[\,\mathrm{Ham}(\sigma(T), \tau(T)) | T = t \text{ and } \neg\mathcal{G}] \le \frac{2\times 0 + 2\times 1 + 2\times 2}{6} = 1,$$

so

$$\begin{aligned}
&\mathbf{E}[\mathrm{Ham}(\sigma(T), \tau(T)) | T = t] \\
&\quad = \Pr(\mathcal{G} | T = t)\mathbf{E}[\,\mathrm{Ham}(\sigma(T), \tau(T)) | T = t \text{ and } \mathcal{G}] \\
&\qquad + \Pr(\neg\mathcal{G} | T = t)\mathbf{E}[\,\mathrm{Ham}(\sigma(T), \tau(T)) | T = t \text{ and } \neg\mathcal{G}] \\
&\quad \le \Pr(\mathcal{G} | T = t)(1 - \tfrac{1}{6}) + (1 - \Pr(\mathcal{G} | T = t)) \\
&\quad = 1 - \tfrac{1}{6}\Pr(\mathcal{G} | T = t).
\end{aligned}$$

Thus if $t \in [\alpha n, n]$,

$$\mathbf{E}[\mathrm{Ham}(\sigma(T), \tau(T)) | T = t] \le 1 - \frac{\delta}{6}.$$

Finally,

$$\begin{aligned}
\mathbf{E}[\mathrm{Ham}(\sigma(T), \tau(T)) | T \le n] &= \sum_{t=1}^{n} \mathbf{E}[\mathrm{Ham}(\sigma(T), \tau(T)) | T = t] \\
&\qquad \times \Pr(T = t | T \le n) \\
&\le 1 - \frac{\delta}{6} \sum_{t=\alpha n}^{n} \Pr(T = t | T \le n) \\
&\le 1 - \frac{\delta}{6} \sum_{t=\alpha n}^{n} \Pr(T = t).
\end{aligned} \tag{19}$$

Since $\alpha \le 1/6$, we have

$$\Pr(T < \alpha n) = 1 - \left(\frac{4n-6}{4n}\right)^{\alpha n} = 1 - \left(1 - \frac{6}{4n}\right)^{\alpha n} \le \frac{6\alpha}{4} \le \frac{1}{4}.$$

Also,

$$\Pr(T > n) = \left(\frac{4n-6}{4n}\right)^{n} = \left(1 - \frac{6}{4n}\right)^{n} \le \exp(-6/4) \le \frac{1}{4}.$$

Thus, (19) gives

$$\mathbf{E}[\,\mathrm{Ham}(\sigma(T),\tau(T))|T \le n] \le 1 - \frac{\delta}{12}. \tag{20}$$

Now Theorem 2.2 of [13] tells us that

$$\mathbf{E}[\mathrm{Ham}(\sigma(n),\tau(n)) - 1] \le \Pr(T \le n)(\mathbf{E}[\mathrm{Ham}(\sigma(T),\tau(T))|T \le n] - 1),$$

and by (20), this is at most $-\delta/12$ so

$$\mathbf{E}[\mathrm{Ham}(\sigma(n),\tau(n))] \le 1 - \frac{\delta}{12}.$$

By the "delayed path coupling lemma" of Czumaj et al. (Lemma 2.1 of [13]), the mixing time satisfies

$$\mathrm{Mix}(\mathcal{M}_{\mathrm{Gl}},\varepsilon) \le \frac{12\log(n\varepsilon^{-1})}{\delta}\,n.$$

In the preceding argument, we assumed that $i \in \{3,\dots,n-2\}$ so that vertices $i-1, i-2$ and $i+1, i+2$ all exist. The argument still goes through if $i$ has fewer neighbors to the left (or right). In that case, we just modify the argument by changing the definition of "good" so that it doesn't mention vertices that don't exist.

Thus, we have proved the following:

THEOREM 12. *Let $G$ be the $n$-vertex path, and let $q = 4$. Consider the Markov chain $\mathcal{M}_{\mathrm{Gl}}$ on the state space $\Omega^+$. Then $\mathrm{Mix}(\mathcal{M}_{\mathrm{Gl}},\varepsilon) \le \frac{12}{\delta} n \log(n\varepsilon^{-1})$, where $\delta$ is the constant mentioned above.*

6.3. *Systematic scan for $q = 4$: $O(\log n)$ sweeps suffice.* We will only define the coupling for pairs $(\sigma,\tau) \in S$. Each such pair disagrees at a single vertex $i$. Thus, when we come to re-color a vertex $j$ during the scan, at most one of $\{j-1, j+1\}$ has a disagreement. The coupling that we will use is as follows. If vertex $j$ is not adjacent to a disagreement, then we use the same colors in both copies. On the other hand, if (say) vertex $j-1$ has a disagreement, then we couple the choice of $\sigma_{j-1}$ for $\sigma_j$ and $\tau_{j-1}$ for $\tau_j$ and we couple the choice of $\tau_{j-1}$ for $\sigma_j$ and $\sigma_{j-1}$ for $\tau_j$. Otherwise, we choose the same color in both copies. The coupling if $j+1$ has a disagreement is similar.

In the following sequence of lemmas, we let $i$ denote the rightmost vertex where there is a disagreement between the colorings $\sigma$ and $\tau$. Lemmas 13 through 16 are valid for any $q \geq 4$, and we state them in terms of $q$ so that we can re-use them later for $q > 4$. The first lemma, Lemma 13, analyzes a scan starting from vertex $i+1$.

LEMMA 13. *Suppose that $\sigma$ and $\tau$ differ at vertex $i < n$ and agree to the right of vertex $i$. Obtain $\sigma'$ and $\tau'$ by scanning left to right, starting at vertex $i+1$. Then*

$$\mathbf{E}[\mathrm{Ham}(\sigma', \tau')] - \mathrm{Ham}(\sigma, \tau) \leq \frac{1}{q-1}.$$

PROOF. If $z = i + \ell$ for $\ell \in \{1, \ldots, n-i\}$, then the probability that vertex $z$ becomes a disagreement after the recoloring is $(1/q)^{\ell}$. Thus, the expected number of additional disagreements is

$$\left(\frac{1}{q}\right)^1 + \left(\frac{1}{q}\right)^2 + \cdots + \left(\frac{1}{q}\right)^{n-i} \leq \frac{1}{q} \times \frac{1}{1-1/q} = \frac{1}{q-1}. \qquad \square$$

The next two lemmas analyze a scan starting from vertex $i$.

LEMMA 14. *Suppose $(\sigma, \tau) \in S$ differ on vertex $i$. Let $C = |\{\sigma_{i-1}, \sigma_{i+1}\}|$. ($C$ is the number of colors that are used at neighbors of $i$ in coloring $\sigma$.) Obtain $\sigma'$ and $\tau'$ by scanning left to right, starting at vertex $i$. Then $\mathbf{E}[\mathrm{Ham}(\sigma', \tau')] \leq C/(q-1)$.*

PROOF. Consider the recoloring of vertex $i$ in copy $\sigma$. With probability $1 - C/q$, the chosen color is not in $\{\sigma_{i-1}, \sigma_{i+1}\}$ so $\mathrm{Ham}(\sigma', \tau') = 0$. On the other hand, no matter what color is chosen for vertex $i$, Lemma 13 guarantees that (conditioned on this choice) $\mathbf{E}[\mathrm{Ham}(\sigma', \tau')] \leq 1 + 1/(q-1)$. Thus, we have

$$\mathbf{E}[\mathrm{Ham}(\sigma', \tau')] \leq \frac{C}{q}\left(1 + \frac{1}{q-1}\right) = \frac{C}{q-1}. \qquad \square$$

LEMMA 15. *Suppose colorings $\sigma$ and $\tau$ differ just on vertices $i-1$ and $i$. Obtain $\sigma'$ and $\tau'$ by scanning left to right, starting at vertex $i$. Then*

$$\mathbf{E}[\mathrm{Ham}(\sigma', \tau')] \leq 1 + \frac{3}{q-1}.$$

PROOF. Consider the recoloring of vertex $i$ in copy $\sigma$. With probability at least $1-3/q$, the chosen color is not in $\{\sigma_{i-1}, \tau_{i-1}, \sigma_{i+1}\}$, so $\mathrm{Ham}(\sigma', \tau') = 1$. On the other hand, no matter what color is chosen for vertex $i$, Lemma 13 guarantees that $\mathbf{E}[\mathrm{Ham}(\sigma', \tau')] \le 2 + 1/(q-1)$. Thus, we have

$$\mathbf{E}[\mathrm{Ham}(\sigma', \tau')] \le \left(1 - \frac{3}{q}\right) \cdot 1 + \frac{3}{q}\left(2 + \frac{1}{q-1}\right),$$

which simplifies to the claimed upper bound. □

The next three lemmas analyze a scan starting from vertex $\max\{1, i-1\}$.

LEMMA 16. *Suppose $(\sigma, \tau) \in S$ differ on vertex $i$. Obtain $\sigma'$ and $\tau'$ by scanning left to right, starting at vertex* $\max\{1, i-1\}$. *Then*

$$\mathbf{E}[\mathrm{Ham}(\sigma', \tau')] \le \frac{3}{q-1}.$$

PROOF. If $i = 1$, then Lemma 14 with $C = 1$ shows $\mathbf{E}[\mathrm{Ham}(\sigma', \tau')] \le 1/(q-1)$, which is at most the expression given in the statement of the lemma. Suppose $i > 1$. Consider the recoloring of vertex $i-1$ in copy $\sigma$. With probability $1/q$, color $\tau_i$ is chosen. By Lemma 15, $\mathbf{E}[\mathrm{Ham}(\sigma', \tau')] \le 1 + 3/(q-1)$. Otherwise, $\sigma'_{i-1} = \tau'_{i-1}$, so Lemma 14 guarantees that $\mathbf{E}[\mathrm{Ham}(\sigma', \tau')] \le 2/(q-1)$. Hence,

$$\mathbf{E}[\mathrm{Ham}(\sigma', \tau')] \le \frac{1}{q}\left(1 + \frac{3}{q-1}\right) + \left(1 - \frac{1}{q}\right)\frac{2}{q-1} = \frac{3}{q-1},$$

as claimed. □

For the rest of this section, we restrict attention to the case $q = 4$, which corresponds to the "break even" situation in Lemma 16.

LEMMA 17. *Suppose $(\sigma, \tau) \in S$ differ on vertex $i < n$. Suppose that $\sigma_{i+1} \notin \{\sigma_i, \tau_i, \sigma_{i-2}\}$. Obtain $\sigma'$ and $\tau'$ by scanning left to right, starting at vertex* $\max\{1, i-1\}$. *Then* $\mathbf{E}[\mathrm{Ham}(\sigma', \tau')] \le \frac{11}{12}$.

PROOF. If $i = 1$, then the lemma follows from Lemma 14 with $C = 1$. Suppose $i > 1$. Consider the recoloring of vertex $i-1$ in copy $\sigma$. With probability $\frac{1}{4}$, color $\sigma_{i+1}$ is chosen. The same color is chosen in copy $\tau$, and Lemma 14 with $C = 1$ guarantees that $\mathbf{E}[\mathrm{Ham}(\sigma', \tau')] \le \frac{1}{3}$. With probability $\frac{1}{2}$, the color chosen for vertex $i-1$ is not in $\{\sigma_{i+1}, \tau_i\}$, so $\sigma'_{i-1} = \tau'_{i-1}$. By Lemma 14 with $C = 2$, $\mathbf{E}[\mathrm{Ham}(\sigma', \tau')] \le \frac{2}{3}$. Otherwise, Lemma 15 guarantees that $\mathbf{E}[\mathrm{Ham}(\sigma', \tau')] \le 2$. Thus, we have $\mathbf{E}[\mathrm{Ham}(\sigma', \tau')] \le \frac{1}{4} \cdot \frac{1}{3} + \frac{1}{2} \cdot \frac{2}{3} + \frac{1}{4} \cdot 2 = \frac{11}{12}$. □

LEMMA 18. *Suppose $(\sigma,\tau)\in S$ differ on vertex $i<n$. Suppose that $\sigma_{i+1}=\sigma_i$. Suppose that $\sigma_i \neq \sigma_{i-2}$ and $\tau_i \neq \sigma_{i-2}$. Obtain $\sigma'$ and $\tau'$ by scanning left to right, starting at vertex $\max\{1,i-1\}$. Then $\mathbf{E}[\mathrm{Ham}(\sigma',\tau')] \leq \frac{11}{12}$.*

PROOF. If $i=1$, then the lemma follows from Lemma 14 with $C=1$.

Suppose $i>1$. Consider the recoloring of vertex $i-1$. With probability $\frac{1}{4}$, color $\tau_i$ is chosen in copy $\sigma$ and $\sigma_i$ is chosen in copy $\tau$. Both of these choices are accepted. In this case, consider the recoloring of vertex $i$. With probability $\frac{1}{2}$, the color chosen is not in $\{\sigma_i,\tau_i\}$ and is accepted in both copies, leaving $\mathrm{Ham}(\sigma',\tau')=1$. Otherwise, by Lemma 13, $\mathbf{E}[\mathrm{Ham}(\sigma',\tau')] \leq \frac{7}{3}$. Thus, conditioned on this color choice for vertex $i-1$, we have $\mathbf{E}[\mathrm{Ham}(\sigma',\tau')] \leq \frac{1}{2}\cdot 1+\frac{1}{2}\cdot\frac{7}{3}=\frac{5}{3}$. For any other choice at vertex $i-1$, Lemma 14 guarantees that $\mathbf{E}[\mathrm{Ham}(\sigma',\tau')] \leq \frac{2}{3}$. We conclude that $\mathbf{E}[\mathrm{Ham}(\sigma',\tau')] \leq \frac{1}{4}\cdot\frac{5}{3}+\frac{3}{4}\cdot\frac{2}{3}=\frac{11}{12}$. □

For the remaining lemmas, we analyze a scan starting from vertex 1. These three lemmas imply the result.

LEMMA 19. *Suppose $(\sigma,\tau)\in S$ differ on vertex $n$. Obtain $\sigma'$ and $\tau'$ by scanning left to right, starting at vertex 1. Then $\mathbf{E}[\mathrm{Ham}(\sigma',\tau')] \leq \frac{11}{16}$.*

PROOF. Consider the recoloring of vertex $n-1$ in coloring $\sigma$. With probability $\frac{1}{4}$, color $\tau_n$ is chosen. In this case, $\mathrm{Ham}(\sigma',\tau')\leq 2$. Otherwise, $\sigma'_{n-1}=\tau'_{n-1}$ so $\mathbf{E}[\mathrm{Ham}(\sigma',\tau')] \leq \frac{1}{4}$. Thus, $\mathbf{E}[\mathrm{Ham}(\sigma',\tau')] \leq \frac{1}{4}\cdot 2+\frac{3}{4}\cdot\frac{1}{4}=\frac{11}{16}$. □

LEMMA 20. *Suppose $(\sigma,\tau)\in S$ differ on vertex $i<n$. Suppose that $\sigma_{i+1}\notin\{\sigma_i,\tau_i\}$. Obtain $\sigma'$ and $\tau'$ by scanning left to right, starting at vertex 1. Then $\mathbf{E}[\mathrm{Ham}(\sigma',\tau')] \leq 47/48$.*

PROOF. If $i\leq 2$, then the lemma follows from Lemma 17. Suppose $i>2$ and consider the recoloring of vertex $i-2$. With probability $\frac{1}{4}$, the color that is chosen is the first color that is not in $\{\sigma'_{i-3},\sigma_{i-1},\sigma_{i+1}\}$. This is accepted so Lemma 17 guarantees that $\mathbf{E}[\mathrm{Ham}(\sigma',\tau')] \leq \frac{11}{12}$. Otherwise, Lemma 16 guarantees that $\mathbf{E}[\mathrm{Ham}(\sigma',\tau')] \leq 1$. Putting this together, $\mathbf{E}[\mathrm{Ham}(\sigma',\tau')] \leq \frac{1}{4}\cdot\frac{11}{12}+\frac{3}{4}\cdot 1=\frac{47}{48}$. □

LEMMA 21. *Suppose $(\sigma,\tau)\in S$ differ on vertex $i<n$. Suppose that $\sigma_{i+1}=\sigma_i$. Obtain $\sigma'$ and $\tau'$ by scanning left to right, starting at vertex 1. Then $\mathbf{E}[\mathrm{Ham}(\sigma',\tau')] \leq 191/192$.*

PROOF. If $i \le 2$, then the lemma follows from Lemma 18. Next suppose $i = 3$. Consider the recoloring of vertex $i-2$. With probability $\frac{1}{4}$, it is recolored with the first color that is not in $\{\sigma_{i-1}, \sigma_i, \tau_i\}$. Now Lemma 18 guarantees that $\mathbf{E}[\mathrm{Ham}(\sigma', \tau')] \le \frac{11}{12}$. Otherwise, Lemma 16 guarantees that $\mathbf{E}[\mathrm{Ham}(\sigma', \tau')] \le 1$. Our conclusion for $i = 3$ is that $\mathbf{E}[\mathrm{Ham}(\sigma', \tau')] \le \frac{1}{4} \cdot \frac{11}{12} + \frac{3}{4} \cdot 1 = \frac{47}{48}$. Finally, suppose $i > 3$. Consider the recoloring of vertex $i-3$. Let $c$ be the first color that is not in $\{\sigma_{i-1}, \sigma_i, \tau_i\}$. With probability $\frac{1}{4}$, the color that is chosen for vertex $i-3$ is the first color that is not in $\{\sigma_{i-4}, \sigma_{i-2}, c\}$. Suppose this happens. Then with probability $\frac{1}{4}$, $c$ is chosen for vertex $i-2$. Then Lemma 18 guarantees that $\mathbf{E}[\mathrm{Ham}(\sigma', \tau')] \le \frac{11}{12}$. Otherwise, Lemma 16 guarantees that $\mathbf{E}[\mathrm{Ham}(\sigma', \tau')] \le 1$. We conclude that $\mathbf{E}[\mathrm{Ham}(\sigma', \tau')] \le \frac{1}{4} \cdot \frac{1}{4} \cdot \frac{11}{12} + \frac{15}{16} \cdot 1 = \frac{191}{192}$. $\square$

Lemmas 19, 20 and 21 imply the following result (by path coupling).

THEOREM 22. *Let $G$ be the $n$-vertex path and let $q = 4$. Consider the Markov chain $\mathcal{M}_{\rightarrow}$ on the state space $\Omega^+$. Then* $\mathrm{Mix}(\mathcal{M}_{\rightarrow}, \varepsilon) \le 192 \log(n\varepsilon^{-1})$.

6.4. *Lower bounds for $q \ge 4$.* In this section we prove that Glauber requires $\Omega(n \log n)$ updates and scan requires $\Omega(\log n)$ sweeps. We use the "disagreement percolation" method of van den Berg [28].

6.4.1. *Calculating the stationary distribution for bounded line segments.* Consider an $s$-edge path (for any $s$). Consider the $q \times q$ "transfer matrix"

$$A = \begin{pmatrix} 0 & 1 & 1 & \cdots & 1 & 1 \\ 1 & 0 & 1 & \cdots & 1 & 1 \\ 1 & 1 & 0 & & 1 & 1 \\ \vdots & \vdots & & \ddots & & \vdots \\ 1 & 1 & 1 & & 0 & 1 \\ 1 & 1 & 1 & \cdots & 1 & 0 \end{pmatrix}.$$

Note that $A^s[i,j]$ is the number of colorings of the path in which the right vertex is colored with color $i$ and the left vertex is colored with color $j$. We will write $e_i$ to denote the row vector with a 1 in column $i$ and zeros elsewhere. Write $f$ to denote the row vector $(1, 1, \ldots, 1)$. Write $v_i$ to denote the row vector with $q-1$ in column $i$ and $-1$ elsewhere. Let $e'_i$, $f'$ and $v'_i$ be the corresponding column vectors. Thus, $e_i A^s e'_j$ is the number of colorings from color $i$ on the right to color $j$ on the left.

Now the (right) eigenvectors of $A$ are $f'$ with eigenvalue $q-1$ and, for every $j$, $v'_j$ with eigenvalue $-1$. Since $e_j = q^{-1} f + q^{-1} v_j$, we have

$$A^s e'_j = q^{-1}(A^s f' + A^s v'_j) = q^{-1}((q-1)^s f' + (-1)^s v'_j),$$

by induction. Thus, if $s$ is even, we have

$$A^s e'_j = q^{-1}((q-1)^s f' + v'_j).$$

So, for $i \neq j$, the number of paths from color $i$ to color $j$ is

$$e_i A^s e'_j = q^{-1}((q-1)^s - 1). \tag{21}$$

Also, the number of paths from color $j$ to color $j$ is

$$\begin{aligned} e_j A^s e'_j &= q^{-1}((q-1)^s + q - 1) \\ &= q^{-1}(q-1)^s(1 + (q-1)^{-(s-1)}). \end{aligned} \tag{22}$$

6.4.2. *Calculating the induced distribution on the color of an internal vertex.* Suppose that $\ell$ and $r$ are positive even integers and let $k = \ell + r$. Consider a path on vertices $1, \ldots, 1+k$. Consider the uniform distribution on colorings in which vertices 1 and $1+k$ are both colored with color $j$.

We wish to bound the probability that vertex $1+\ell$ is colored with color $j$. From (22), this is

$$\begin{aligned} &\frac{q^{-1}(q-1)^r(1+(q-1)^{-(r-1)}) \times q^{-1}(q-1)^\ell(1+(q-1)^{-(\ell-1)})}{q^{-1}(q-1)^k(1+(q-1)^{-(k-1)})} \\ &= q^{-1}(1+(q-1)^{-(r-1)}) \frac{1+(q-1)^{-(\ell-1)}}{1+(q-1)^{-(k-1)}} \\ &\geq q^{-1}(1+(q-1)^{-(r-1)}). \end{aligned} \tag{23}$$

6.4.3. *Dividing the line into segments.* Let $r$ be the largest even number not exceeding $\frac{1}{3}\log_{q-1} n$. Let $\ell$ be the smallest even number that is at least $48 \log_e n$. Let $k = \ell + r$ and $m = \lfloor (n-1)/k \rfloor$. For $i \in \{0, \ldots, m-1\}$, let $L_i$ be the vertex $1 + ik$ and $M_i$ be the vertex $1 + ik + \ell$. Finally, let $L_m = 1 + mk$ and $R_i$ be $L_{i+1}$. The idea is to divide the line into line segments. Segment $i$ has left endpoint $L_i$ and "middle" point $M_i$ (which is not quite in the middle!) and right endpoint $R_i$.

Let $Z_i$ be the indicator for the event that vertex $M_i$ is colored with color 0. Let $Z = \sum_{i=0}^{m-1} Z_i$. Of course, the expectations of $Z_i$ and $Z$ are only well defined if we focus attention on a particular distribution over $\Omega$. We will use the notation $\mathbf{E}_\pi(Z_i)$ to refer to the expectation of $Z_i$ in distribution $\pi$.

6.4.4. *Calculating the distribution of $Z$ in stationarity.* Let $\pi$ be the uniform distribution on $\Omega$. We can sample from $\pi$ by filling in the colors from left to right. There are $q(q-1)^{n-1}$ possible colorings in $\Omega$. Given the colors of vertices $1, \ldots, n-v$ (for $v < n$), there are $(q-1)^v$ ways to finish the coloring, and these are chosen uniformly. The probability that $Z_1 = 1$ is $1/q$. For any $i > 1$, we can use (22) [observing that it is (22) that determines the upper bound, and not (21)] to see that the probability that $Z_i = 1$, conditioned on colors $\sigma_1, \dots, \sigma_{M_{i-1}}$, is at most

$$\frac{q^{-1}(q-1)^k(1+(q-1)^{-(k-1)})}{(q-1)^k} = q^{-1}(1+(q-1)^{-(k-1)}).$$

Thus, $Z$ is dominated from above by the sum of $m$ independent Bernoulli random variables with success probability $p = q^{-1}(1+(q-1)^{-(k-1)})$.

Let $\varepsilon = m^{-3/8}$. By a Chernoff bound,

$$\Pr_\pi(Z \geq (1+\varepsilon)mp) \leq \exp(-\varepsilon^2 mp/3). \tag{24}$$

Note that $(q-1)^{k-1} = \omega(n^{1/3})$ and that $m = \Theta(n/\log n)$. Thus, (24) implies

$$\begin{aligned}\Pr_\pi(Z \geq q^{-1}m + \tfrac{1}{2}mn^{-1/3}) &\leq \Pr_\pi(Z \geq (1+\varepsilon)mp) \\ &\leq \exp(-\varepsilon^2 mp/3) = o(1).\end{aligned} \tag{25}$$

6.4.5. *An initial distribution for the Markov chains.* Equation (25) shows that, in the stationary distribution $\pi$, $Z$ is unlikely to exceed $q^{-1}m + \frac{1}{2}mn^{-1/3}$. In this section we will define an initial distribution $\pi_0$ on colorings in $\Omega$.

The idea will be to show that, if $\sigma(0)$ is chosen from $\pi_0$ and $\sigma(0), \sigma(1), \dots, \sigma(t)$ evolves according to the dynamics (either Glauber or scan) and $t$ is too small, then, in the distribution of $\sigma(t)$, $Z$ is likely to exceed $q^{-1}m + \frac{1}{2}mn^{-1/3}$. This allows us to conclude that the chain does not mix by step $t$.

Let $\pi_0$ be the uniform distribution on colorings in which vertices $L_0, \dots, L_m$ are colored 0. By (23), $\mathbf{E}_{\pi_0} Z \geq mq^{-1}(1+(q-1)^{-(r-1)}) \geq mq^{-1}(1+(q-1)n^{-1/3})$. By a Chernoff bound,

$$\Pr_{\pi_0}(Z \geq q^{-1}m + \tfrac{1}{2}mn^{-1/3}) \geq 1 - o(1). \tag{26}$$

6.4.6. *The $t$-step distribution for systematic scan.* Suppose $\sigma(0)$ is chosen from $\pi_0$. Let $\sigma(0), \sigma(1), \dots$ evolve according to the dynamics of $\mathcal{M}_\rightarrow$. Let $\mathcal{L}_{\rightarrow,t}(Z)$ denote the distribution of the random variable $Z$ in the coloring $\sigma(t)$.

Suppose $\tau(0)$ is chosen from $\pi_0$. Let $\tau(0), \tau(1), \dots$ evolve according to the "clamped dynamics" $\mathcal{M}^{\mathrm{c}}_\rightarrow$, which is the same as $\mathcal{M}_\rightarrow$ except that all moves involving vertices $\{L_0, \dots, L_m\}$ are rejected (so the color of these vertices cannot change). Let $\mathcal{L}^{\mathrm{c}}_{\rightarrow,t}(Z)$ denote the distribution of the random variable $Z$ in the coloring $\tau(t)$. By construction, the distribution $\mathcal{L}^{\mathrm{c}}_{\rightarrow,t}(Z)$ is same as the distribution of $Z$ in $\pi_0$. [This follows because $\pi_0$ is the stationary distribution of $\mathcal{M}^{\mathrm{c}}_\rightarrow$. This can be proved as follows, where $\Omega_0$ is the set of colorings in which vertices $L_0, \dots, L_m$ are colored 0. Let $P^{\mathrm{c}}_\rightarrow$ be the transition matrix of $\mathcal{M}^{\mathrm{c}}_{\rightarrow}$ and let $P^{\mathrm{c}}_{\leftarrow}$ be the transition matrix of the reversal. Then any stationary distribution $\pi'$ of $\mathcal{M}^{\mathrm{c}}_{\rightarrow}$ satisfies

$$\pi'(\sigma') = \sum_{\sigma\in\Omega_0} \pi'(\sigma)P^{\mathrm{c}}_{\rightarrow}(\sigma,\sigma') = \sum_{\sigma\in\Omega_0} \pi'(\sigma)P^{\mathrm{c}}_{\leftarrow}(\sigma',\sigma),$$

but the latter equation is satisfied by the uniform distribution $\pi' = \pi_0$. Also, the chain is ergodic so has a unique stationary distribution.]

To upper bound $\mathrm{d_{TV}}(\mathcal{L}_{\rightarrow,t}(Z), \mathcal{L}^{\mathrm{c}}_{\rightarrow,t}(Z))$, we will consider a joint process $(\sigma(t), \tau(t))$ in which the first component has the same distribution as $(\sigma(t))$ and the second component has the same distribution as $(\tau(t))$. The total variation distance $\mathrm{d_{TV}}(\mathcal{L}_{\rightarrow,t}(Z), \mathcal{L}^{\mathrm{c}}_{\rightarrow,t}(Z))$ is upper-bounded by the probability that some vertex $M_i$ gets different colors in $\sigma(t)$ and $\tau(t)$.

The particular joint process that we will consider starts with $\sigma(0) = \tau(0)$. To move from $(\sigma(t-1), \tau(t-1))$ to $(\sigma(t), \tau(t))$, we use the "switch coupling." When we consider vertex $v$ for recoloring, we will couple the color choices as follows:

(A) if we consider color $\sigma(v-1)$ for $v$ in $\sigma$, then consider color $\tau(v-1)$ for $v$ in $\tau$,
(B) if we consider color $\tau(v-1)$ for $v$ in $\sigma$, then consider color $\sigma(v-1)$ for $v$ in $\tau$,
(C) otherwise consider the same color for $v$ in $\tau$ as in $\sigma$.

We will be particularly interested in $t < r$. For such a $t$, and for any $i \in \{0,\dots,m-1\}$, the probability that vertex $M_i$ gets different colors in $\sigma(t)$ and $\tau(t)$ is at most the probability that we chose option (B) in order for vertices $L_i+1,\dots,L_i+\ell$ over the $t$ scans.

Say that vertex $L_i+v$ is "interrupting" (i.e., it interrupts the disagreement percolation) if, the first time that we consider this vertex when we have a disagreement at vertex $L_i+v-1$, we choose some option other than option (B) for vertex $L_i+v$.

The probability that vertex $M_i$ gets different colors in $\sigma(t)$ and $\tau(t)$ is at most the probability that we have fewer than $t$ interrupting vertices in $L_i+1,\dots,L_i+\ell$; this probability is dominated (from above) by the probability of having fewer than $t$ successes in $\ell$ Bernoulli trials with success probability $(q-1)/q$. So if we take any $t \le r/2 \le (2/3)\ell(q-1)/q$, a Chernoff bound says that the probability of having fewer than $t$ interrupting vertices is at most

$$\exp(-(1/3)^2\ell(q-1)/(2q)),$$

which is at most $n^{-2}$ by the definition of $\ell$. Thus, the probability that there exists an $i$ such that vertex $M_i$ gets different colors in $\sigma(t)$ and $\tau(t)$ is at most $mn^{-2} = o(1)$.

Thus, for any $t \le r/2$, $\mathrm{d_{TV}}(\mathcal{L}_{\rightarrow,t}(Z), \mathcal{L}^{\mathrm{c}}_{\rightarrow,t}(Z)) = o(1)$, so, by (26),

$$\Pr_{\mathcal{L}_{\rightarrow,t}(Z)}(Z \ge q^{-1}m + \tfrac{1}{2}mn^{-1/3}) \ge 1 - o(1).$$

Combining this with (25), we find that $d_{TV}(\mathcal{L}_{\to,t},\pi) \geq 1-o(1)$ so systematic scan does not mix in $t$ steps. Thus, we obtain the following theorem.

THEOREM 23. *Let $G$ be the $n$-vertex path, and let $q \geq 4$. Consider the Markov chain $\mathcal{M}_{\to}$ on the state space $\Omega$. For any fixed $\varepsilon < 1$ and sufficiently large $n$,*

$$\mathrm{Mix}(\mathcal{M}_{\to},\varepsilon) \geq \tfrac{1}{2}(\tfrac{1}{3}\log_{q-1} n - 2).$$

6.4.7. *The $t$-step distribution for Glauber.* A similar argument to that of Section 6.4.6 can be used to show that Glauber dynamics does not mix in $t$ steps for some $t = \Omega(n \log_{q-1} n)$. The particular value of $t$ for which the straightforward argument works is around $nr/(q^2 e)$. We prefer to give a stronger argument which gives a better bound as a function of $q$. The idea for the stronger argument is as follows. In Section 6.4.6 we showed that the distribution of $\mathcal{L}_{\to,t}(Z)$ and $\mathcal{L}^{c}_{\to,t}(Z)$ were close by showing that, with high probability, there was no $i \in \{0,\dots,m-1\}$ for which a disagreement at vertex $L_i$ or $R_i$ could percolate to vertex $M_i$. Here we observe that the distributions $\mathcal{L}_{\to,t}(Z)$ and $\mathcal{L}^{c}_{\to,t}(Z)$ would be close even if some, but not many, of the percolations occur.

We start with some notation. It will be helpful to keep track of the nearest endpoint to an arbitrary vertex $v$. For this purpose, if $v$ is in the range $L_i+1,\dots,M_i-1$, its "*important neighbor*" will be vertex $v-1$. On the other hand, if $v$ is in the range $M_i,\dots,R_i-1$, its important neighbor will be vertex $v+1$.

As in Section 6.4.6, we will consider a process $\sigma(0),\sigma(1),\dots$ in which $\sigma(0)$ is drawn from $\pi_0$ and $\sigma(t)$ evolves according to $\mathcal{M}_{\mathrm{Gl}}$. We will also consider the process $\tau(0),\tau(1),\dots$ in which $\tau(0)=\sigma(0)$ and $\tau(t)$ evolves according to a *clamped dynamics* $\mathcal{M}^{c}_{\mathrm{Gl}}$ in which moves involving $L_0,L_1,\dots,L_m$ are rejected. We will construct a joint process $(\sigma(t),\tau(t))$ with $\sigma(0)=\tau(0)$. To move from $(\sigma(t-1),\tau(t-1))$ to $(\sigma(t),\tau(t))$, we choose the same vertex $v$ in both copies. If $v=L_i$, for some $i$ then only $\sigma$ is changed. If $v>L_m$, then we use the same color in both copies. Otherwise, we do a switch coupling, based on the important neighbor, $w$, of $v$. In particular, we couple the color choices as follows:

(A) if we consider color $\sigma(w)$ for $v$ in $\sigma$, then consider color $\tau(w)$ for $v$ in $\tau$,
(B) if we consider color $\tau(w)$ for $v$ in $\sigma$, then consider color $\sigma(w)$ for $v$ in $\tau$,
(C) otherwise consider the same color for $v$ in $\tau$ as in $\sigma$.

Suppose that $t \leq (qnr)/(2e(q-1))$. The probability that $M_i$ gets different colors in $\sigma(t)$ and $\tau(t)$ is at most the probability that (at least) one of the following occurs during the $t$ steps:

- During some ordered sequence of $\ell-1$ steps, the process recolors vertices $L_i+1,\ L_i+2,\dots,L_i+\ell-1=M_i-1$ using option (B).
- During some ordered sequence of $r$ steps, the process recolors vertices $R_i-1,\ R_i-2,\dots,R_i-r=M_i$ using option (B).

The probability that one of these occurs is at most

$$2\binom{t}{r}\left(\frac{1}{qn}\right)^r \le 2\left(\frac{te}{rqn}\right)^r \le \frac{1}{8}\left(\frac{2te}{rqn}\right)^r \le \frac{1}{8}\left(\frac{1}{q-1}\right)^r,$$

where we have crudely used $r\ge 4$ in the second inequality. Thus,

$$\mathrm{d_{TV}}(\mathcal{L}_{\mathrm{Gl},t}(Z_i),\mathcal{L}^{\mathrm{c}}_{\mathrm{Gl},t}(Z_i)) \le \frac{1}{8}\left(\frac{1}{q-1}\right)^r. \tag{27}$$

Now combining (23) and (27) we have

$$\begin{aligned}\Pr_{\mathcal{L}_{\mathrm{Gl},t}}(Z_i=1) &\ge \Pr_{\mathcal{L}^{\mathrm{c}}_{\mathrm{Gl},t}}(Z_i=1) - \mathrm{d_{TV}}(\mathcal{L}_{\mathrm{Gl},t}(Z_i),\mathcal{L}^{\mathrm{c}}_{\mathrm{Gl},t}(Z_i))\\ &\ge q^{-1}(1+(q-1)^{-(r-1)}) - \frac{1}{8}\left(\frac{1}{q-1}\right)^r\\ &\ge q^{-1}+\frac{5}{8}(q-1)^{-r}\\ &\ge q^{-1}+\frac{5}{8}n^{-1/3}.\end{aligned}$$

So

$$\mathbf{E}_{\mathcal{L}_{\mathrm{Gl},t}}(Z) \ge q^{-1}m+\tfrac{5}{8}mn^{-1/3}.$$

Also,

$$\mathrm{var}_{\mathcal{L}_{\mathrm{Gl},t}}(Z_i) = \Pr_{\mathcal{L}_{\mathrm{Gl},t}}(Z_i=1)\Pr_{\mathcal{L}_{\mathrm{Gl},t}}(Z_i=0)\le 1.$$

We will show in Lemma 25 (below) that, for $i\ne j$, $\mathrm{cov}_{\mathcal{L}_{\mathrm{Gl},t}}(Z_i,Z_j)\le m^{-1}$. So

$$\begin{aligned}\mathrm{var}_{\mathcal{L}_{\mathrm{Gl},t}}(Z) &= \sum_i \mathrm{var}_{\mathcal{L}_{\mathrm{Gl},t}}(Z_i) + \sum_{i\ne j}\mathrm{cov}_{\mathcal{L}_{\mathrm{Gl},t}}(Z_i,Z_j)\\ &\le m + \sum_{i\ne j}\mathrm{cov}_{\mathcal{L}_{\mathrm{Gl},t}}(Z_i,Z_j)\\ &\le 2m.\end{aligned}$$

Let

$$\lambda = \frac{(1/8)mn^{-1/3}}{\sqrt{2m}}.$$

Note that $\lambda = \omega(1)$ as a function of $n$. Also,

$$\mathbf{E}_{\mathcal{L}_{\mathrm{Gl},t}}(Z) - \lambda\sqrt{\mathrm{var}_{\mathcal{L}_{\mathrm{Gl},t}}(Z)} \geq q^{-1}m + \tfrac{5}{8}mn^{-1/3} - \lambda\sqrt{2m}$$
$$= q^{-1}m + \tfrac{1}{2}mn^{-1/3}.$$

Thus, by Chebyshev's inequality, we have

$$\Pr_{\mathcal{L}_{\mathrm{Gl},t}}(Z \leq q^{-1}m + \tfrac{1}{2}mn^{-1/3}) \leq \Pr_{\mathcal{L}_{\mathrm{Gl},t}}\left(Z \leq \mathbf{E}_{\mathcal{L}_{\mathrm{Gl},t}}(Z) - \lambda\sqrt{\mathrm{var}_{\mathcal{L}_{\mathrm{Gl},t}}(Z)}\right)$$
$$\leq \lambda^{-2} = o(1).$$

Combining this with (25), we find that $d_{\mathrm{TV}}(\mathcal{L}_{\mathrm{Gl},t}, \pi) \geq 1 - o(1)$, so Glauber dynamics does not mix in $t$ steps for any $t \leq (qnr)/(2e(q-1))$. Thus, we obtain the following theorem.

THEOREM 24. *Let $G$ be the $n$-vertex path, and let $q \geq 4$. Consider the Markov chain $\mathcal{M}_{\mathrm{Gl}}$ on the state space $\Omega$. For any fixed $\varepsilon < 1$ and sufficiently large $n$,*

$$\mathrm{Mix}(\mathcal{M}_{\mathrm{Gl}}, \varepsilon) \geq \frac{qn((1/3)\log_{q-1} n - 2)}{2e(q-1)}.$$

LEMMA 25. *For $i \neq j$, $\mathrm{cov}_{\mathcal{L}_{\mathrm{Gl},t}}(Z_i, Z_j) \leq m^{-1}$.*

PROOF. We will show that $Z_i$ and $Z_j$ have low covariance in the $t$-step distribution by showing that Glauber dynamics (over $t$ steps) is quite close to a "clamped distribution" in which some vertex between $M_i$ and $M_j$ is held fixed. This "disagreement percolation" argument is similar to the argument in Section 6.4.7. The only difference is that, in order to get a sufficiently small upper bound on the covariance, we have to look at a "clamped process" that is slightly different from $\mathcal{M}_{\mathrm{Gl}}^{\mathrm{c}}$. In particular, in $\mathcal{M}_{\mathrm{Gl}}^{\mathrm{c}}$, $M_i$ is only $r$ vertices away from the nearest "clamped vertex," $R_i$. Here we need to spread the clamped vertices out more symmetrically with respect to the vertices $M_i$. Let

$$\Gamma = \{M_i + k/2 | i \in \{0, \dots, m-1\}\}.$$

Consider a process $\sigma(0), \sigma(1), \dots$ which evolves according to $\mathcal{M}_{\mathrm{Gl}}$. Let $\rho(0)$, $\rho(1), \dots$ be a process which evolves according to a clamped version of $\mathcal{M}_{\mathrm{Gl}}$ in which those moves involving vertices in $\Gamma$ are rejected. We refer to this dynamics as $\mathcal{M}_{\mathrm{Gl}}^{\mathrm{sc}}$ (where "sc" is intended to indicate "symmetric clamped"). Consider the joint process $(\sigma(t), \rho(t))$ which starts with $\rho(0) = \sigma(0)$ and progresses according to the identity coupling [the same vertices and colors are chosen in the transition $\sigma(t-1) \to \sigma(t)$ and in the transition $\rho(t-1) \to$ $\rho(t)]$. Now the probability that $\sigma(t)_{M_i} \neq \rho(t)_{M_i}$ or $\sigma(t)_{M_j} \neq \rho(t)_{M_j}$ (or both) is at most

$$4\binom{t}{k/2}\left(\frac{1}{n}\right)^{k/2},$$

since an ordered sequence of $k/2$ vertices would need to be chosen either from the left toward $M_i$ or from the right toward $M_i$ or from the left or right toward $M_j$. The probability that a particular vertex is chosen at any step is $1/n$. Since $t \leq (qnr)/(2e(q-1))$ and $k \geq 8er$, this is at most

$$\left(\frac{8te}{kn}\right)^{k/2} \leq e^{-k/2} \leq n^{-24}.$$

Now

$$\begin{aligned}
\operatorname{cov}_{\mathcal{L}_{\mathrm{Gl},t}}(Z_i, Z_j) &= \mathbf{E}_{\mathcal{L}_{\mathrm{Gl},t}}(Z_i Z_j) - \mathbf{E}_{\mathcal{L}_{\mathrm{Gl},t}}(Z_i)\mathbf{E}_{\mathcal{L}_{\mathrm{Gl},t}}(Z_j) \\
&= \Pr_{\mathcal{L}_{\mathrm{Gl},t}}(Z_i = 1 \wedge Z_j = 1) - \Pr_{\mathcal{L}_{\mathrm{Gl},t}}(Z_i = 1)\Pr_{\mathcal{L}_{\mathrm{Gl},t}}(Z_j = 1) \\
&\leq \Pr_{\mathcal{L}^{\mathrm{sc}}_{\mathrm{Gl},t}}(Z_i = 1 \wedge Z_j = 1) + n^{-24} \\
&\quad - (\Pr_{\mathcal{L}^{\mathrm{sc}}_{\mathrm{Gl},t}}(Z_i = 1) - n^{-24})(\Pr_{\mathcal{L}^{\mathrm{sc}}_{\mathrm{Gl},t}}(Z_j = 1) - n^{-24}) \\
&\leq \operatorname{cov}_{\mathcal{L}^{\mathrm{sc}}_{\mathrm{Gl},t}}(Z_i, Z_j) + 4n^{-24} \\
&= 4n^{-24}.
\end{aligned}$$

□

**7. Optimal mixing for Glauber and scan when $q > 4$.** Let $G$ be the $n$-vertex path. For $q > 4$, Lemma 1 of [22] shows that Glauber dynamics mixes in $O(n \log n)$ steps. For scan, we use the coupling from Section 6.3. Consider a pair $(\sigma, \tau) \in S$ which disagrees at a single vertex $i$. Obtain $\sigma'$ and $\tau'$ by scanning left to right, starting at vertex $\max\{1, i-1\}$. Lemma 16 shows that $\mathbf{E}[\mathrm{Ham}(\sigma', \tau')] \leq 3/4$. This implies the following theorem (by path coupling).

THEOREM 26. *Let $G$ be the $n$-vertex path, and let $q > 4$. Consider the Markov chain $\mathcal{M}_{\rightarrow}$ on the state space $\Omega^+$. Then* $\mathrm{Mix}(\mathcal{M}_{\rightarrow}, \varepsilon) \leq 4\log(n\varepsilon^{-1})$.

**8. $H$-coloring: $O(n^5)$ updates or scans suffice.** Let $H$ be a fixed graph, possibly with self-loops, and let $\Omega$ be the set of $H$-colorings of the graph $G$. These are the homomorphisms from $G$ to $H$—see [7, 15, 19] for details. We can extend the dynamics $\mathcal{M}_{\mathrm{Gl}}$ and $\mathcal{M}_{\rightarrow}$ to the domain of $H$-coloring by modifying the procedure Metropolis$(v)$ from Section 2. In particular, a proposed color $c$ (which is a vertex of $H$) is accepted if and only if every neighbor $w$ of $v$ is colored with some neighbor $c_w$ of $c$. The original dynamics corresponds to the situation in which $H$ is a $q$-clique with no self-loops.

Suppose that $H$ is connected. Let $G$ be the $n$-vertex path. If $H$ has an odd cycle, then Glauber dynamics and systematic scan are both ergodic on $\Omega$, the set of $H$-colorings of $G$. In this case we say that any two colorings, $\sigma \in \Omega$ and $\tau \in \Omega$, are *compatible*. If $H$ does not have an odd cycle, then it is bipartite. Neither dynamics is ergodic on $\Omega$. However, the $H$-colorings can be partitioned in a natural way into two subsets, such that Glauber and scan are both ergodic on either subset. In particular, the $H$-colorings are partitioned as follows. Two $H$-colorings $\sigma \in \Omega$ and $\tau \in \Omega$ are compatible if $\sigma_1$ and $\tau_1$ are chosen from the same side of the bipartition of $H$. Our aim is to show rapid mixing on the set(s) of compatible colorings:

Let $h = |V(H)|$. Define $t$ as follows:

$$t = \begin{cases} 4h-1, & \text{if } H \text{ is not bipartite and } n \text{ is even;} \\ 2h-1, & \text{if } H \text{ is bipartite and } n \text{ is even;} \\ 4h, & \text{if } H \text{ is not bipartite and } n \text{ is odd;} \\ 2h, & \text{if } H \text{ is bipartite and } n \text{ is odd.} \end{cases}$$

Note that $n+t$ is always odd.

LEMMA 27. *In any two compatible $H$-colorings $\sigma$ and $\tau$, there is a $t$-edge path in $H$ from $\sigma_n$ to $\tau_1$.*

PROOF. We look at each of the four cases.

$H$ is not bipartite and $n$ is even: Let $c$ be some point on an odd-length cycle. Go from $\sigma_n$ to $c$ in at most $h-1$ edges. Also, go from $c$ to $\tau_1$ in at most $h-1$ edges. If the constructed path has an even number of edges, go around the cycle using at most $h$ more edges. Now go back and forth on the last edge to make the total length equal to $t$.

$H$ is bipartite and $n$ is even: Note that $\sigma_n$ and $\tau_1$ are on opposite sides of the bipartition. Go from $\sigma_n$ to $\tau_1$ in at most $h-1$ edges and go back and forth on the last edge.

$H$ is not bipartite and $n$ is even: Let $c$ be some point on an odd-length cycle. Go from $\sigma_n$ to $c$ in at most $h-1$ edges. Also, go from $c$ to $\tau_1$ in at most $h-1$ edges. If the constructed path has an odd number of edges, go around the cycle using at most $h$ more edges. Now go back and forth on the last edge to make the total length equal to $t$.

$H$ is bipartite and $n$ is odd: Note that $\sigma_n$ and $\tau_1$ are on the same side of the bipartition. Go from $\sigma_n$ to $\tau_1$ in at most $h-1$ edges and go back and forth on the last edge. $\square$

8.1. *Constructing canonical paths.* Let $\Omega'$ be the state space. It is either the set of all proper colorings (if $H$ is not bipartite) or it is one of the two maximum sets of compatible colorings (if $H$ is bipartite). We will use the canonical paths method, which can be viewed as a special case of comparison in which we compare $\mathcal{M}_{\mathrm{Gl}}$ to the uniform random walk on $\Omega'$. Thus, for each $\sigma \in \Omega'$ and $\tau \in \Omega'$, we will construct a canonical path $\gamma_{\sigma,\tau}$ from $\sigma$ to $\tau$.

First, let $\sigma_n c_1 \cdots c_{t-1} \tau_1$ be some $t$-edge path from $\sigma_n$ to $\tau_1$ and let $z_1 z_2 \cdots z_{2n+t-1}$ denote $\sigma_1 \cdots \sigma_n c_1 \cdots c_{t-1} \tau_1 \cdots \tau_n$. Let $Z_i$ denote the $H$-coloring $z_i z_{i+1} \cdots z_{i+n-1}$, so $Z_1 = \sigma$ and $Z_{n+t} = \tau$. The path $\gamma_{\sigma,\tau}$ passes through $Z_1, Z_3, Z_5, \ldots, Z_{n+t}$. Moving from $Z_i$ to $Z_{i+2}$ can be implemented by $n$ Glauber transitions (applied to vertices 1 to $n$ in order). Let

$$A = \max_{\alpha,\beta} \frac{1}{\pi(\alpha) P_{\mathrm{Gl}}(\alpha,\beta)} \sum_{\sigma,\tau} \pi(\sigma)\pi(\tau)\,|\gamma_{\sigma,\tau}|, \tag{28}$$

where the max is over all Glauber-dynamics transitions $(\alpha, \beta)$ and the sum is over all pairs $(\sigma, \tau)$ such that $(\alpha, \beta)$ is on the canonical path $\gamma_{\sigma,\tau}$. By Theorem 2.1 of [9], we have $\lambda(\mathcal{M}_{\mathrm{Gl}}) \geq 1/A$. We now derive an upper bound on $A$.

The three stationary probabilities in (28) are all $1/|\Omega'|$. Furthermore, every canonical path $\gamma_{\sigma,\tau}$ satisfies $|\gamma_{\sigma,\tau}| \leq \frac{n+t}{2} n$. Finally, $P_{\mathrm{Gl}}(\alpha, \beta) = \frac{1}{nh}$. Plugging this into (28), we get

$$A \leq \frac{n+t}{2} n \frac{nh}{|\Omega'|} \max_{\alpha,\beta} \sum_{\sigma,\tau} 1.$$

We will show that the number of pairs $(\sigma, \tau)$ using transition $(\alpha, \beta)$ is $O(n\,|\Omega'|)$, from which we can conclude

$$A = O(n^4), \tag{29}$$

viewing $h$ as constant. The method we use is standard: each canonical path through $(\alpha, \beta)$ will be assigned a unique "encoding" chosen from a set of $O(n\,|\Omega'|)$ encodings.

So now fix $(\alpha, \beta)$ and consider the set of all canonical paths that use transition $(\alpha, \beta)$. We show how to encode a typical such path, from $\sigma$ to $\tau$, say. Let $\tau_n c'_1 \cdots c'_{t-1} \sigma_1$ be some $t$-edge path in $H$ from $\tau_n$ to $\sigma_1$ and let $\hat{z}_1 \hat{z}_2 \cdots \hat{z}_{2n+t-1}$ denote the path $\tau_1 \cdots \tau_n c'_1 \cdots c'_{t-1} \sigma_1 \cdots \sigma_n$. Let $\hat{Z}_i$ denote the $H$-coloring $\hat{z}_i \hat{z}_{i+1} \cdots \hat{z}_{i+n-1}$.

The encoding of the canonical path from $\sigma$ to $\tau$ consists of the following information:

- $i$, indicating that the current transition is on the path from $Z_i$ to $Z_{i+2}$, and the colors $z_i$ and $z_{i+1}$;
- $\hat{Z}_i$, and
- the colors $\sigma_{i-t+1}, \ldots, \sigma_{i-1}$ and $\tau_{i-t+1}, \ldots, \tau_{i-1}$.

Given the transition $(\alpha, \beta)$ and the values of $i$, $z_i$ and $z_{i+1}$, we can deduce $Z_i$. From $Z_i$ and $\hat{Z}_i$ and the extra colors, we can deduce $\sigma$ and $\tau$. Thus, the number of pairs $(\sigma, \tau)$ using the given transition is at most the number of encodings, which is $O(n\,|\Omega'|)$ as required, so we have now established (29).

Note that $\mathcal{M}_{\mathrm{Gl}}$ is reversible. Let $1 = \beta_0 > \beta_1 \geq \cdots \geq \beta_{|\Omega'|-1} > -1$ be the eigenvalues of its transition matrix $P_{\mathrm{Gl}}$. Since $1 - \beta_1 = \lambda(\mathcal{M}_{\mathrm{Gl}})$, we have $1/(1-\beta_1) \leq A$. To bound the mixing time of $\mathcal{M}_{\mathrm{Gl}}$, we also need an upper bound on $\frac{1}{1+\beta_{|\Omega'|-1}}$. This is an easy application of Proposition 2 of [10] since, for every $\sigma \in \Omega'$, we have $P_{\mathrm{Gl}}(\sigma,\sigma) \geq 1/h$.

In particular, for every $\sigma \in \Omega'$, we define the (odd-length) canonical path from $\sigma$ to itself to be single transition $P_{\mathrm{Gl}}(\sigma,\sigma)$. Proposition 2 of [10] then gives

$$\frac{1}{1+\beta_{|\Omega'|-1}} \leq \frac{1}{2} \max_{\sigma} \frac{1}{P_{\mathrm{Gl}}(\sigma,\sigma)} \leq \frac{h}{2}.$$

Combining this with (29), Proposition 1(i) of [26] gives

$$\begin{aligned}\mathrm{Mix}_\sigma(\mathcal{M}_{\mathrm{Gl}}, \varepsilon) &\leq \frac{1}{1-\beta_{\max}} (\ln \pi(\sigma)^{-1} + \ln \varepsilon^{-1}) \\ &= O(n^4)(\ln \pi(\sigma)^{-1} + \ln \varepsilon^{-1}).\end{aligned}$$

Thus, we have the following theorem.

Theorem 28. *Let $H$ be a fixed connected graph. Let $G$ be the $n$-vertex path. Let $\Omega'$ be the state space of $\mathcal{M}_{\mathrm{Gl}}$, which is either the set of all proper $H$-colorings of $G$ (if $H$ is not bipartite) or one of the two maximum sets of compatible colorings (if $H$ is bipartite). Consider the Markov chain $\mathcal{M}_{\mathrm{Gl}}$ on the state space $\Omega'$. Then*

$$\mathrm{Mix}(\mathcal{M}_{\mathrm{Gl}}, \varepsilon) = O(n^5 \ln \varepsilon^{-1}).$$

In Section 10.1 we will show how to use our lower bound $\lambda(\mathcal{M}_{\mathrm{Gl}}) \geq 1/A$ to get a corresponding lower bound on $\lambda(\mathcal{M}_{\rightarrow})$. This will imply that the mixing time of systematic scan is also $O(n^5 \ln \varepsilon^{-1})$, though, for technical reasons (since scan is not reversible), we state the result in continuous time. See Theorem 31 in Section 10.1 for details.

8.2. *Special case.* Suppose that $H$ is an odd cycle of length $k$. We noted at the beginning of Section 8 that Glauber and scan are ergodic on $\Omega$ and Section 8.1 shows that the mixing time is $O(n^5)$. In fact, the analysis for 3-coloring translates directly to the case of a $k$-cycle so we get the following analog of Theorems 1 and 3.

Theorem 29. *Let $H$ be an odd cycle. Let $G$ be the $n$-vertex path. Let $\Omega'$ be the set of $H$-colorings of $G$. Consider the Markov chain $\mathcal{M}_{\mathrm{Gl}}$ on the state space $\Omega'$.* $\mathrm{Mix}(\mathcal{M}_{\mathrm{Gl}}, \frac{1}{2}) = \Theta(n^3 \log n)$.

The generalization of the proofs of Theorems 1 and 3 is straightforward. In Section 4.1 each configuration $X \in \Upsilon$ corresponds to $k$ colorings. In Section 4.2 the height $h_i$ of every vertex $i$ satisfies $h_i = \sigma_i \pmod k$. The quantity $B$ in Section 4.3 is increased by a factor of $k$. A similar result holds for scan.

**9. Directed *H*-coloring.** It is natural to ask whether the $H$-coloring results could be generalized, for example, to directed $H$-coloring. The answer is no. To illustrate this, we give an example of a directed $H$ that is not ergodic on the $n$-vertex path $G$, and another example of a directed $H$ for which Glauber is ergodic, but mixes slowly.

For the first example, let $H$ have vertex set $\{x, y, z\}$ and edge set $\{(x, y), (y, z), (z, x)\}$. Now the three possible colorings of $G$ are

$$xyzxyz\ldots, \qquad yzxyzx\ldots \quad \text{and} \quad zxyzxy\ldots.$$

These are not connected by either Glauber or scan moves.

For the second example, let the vertices of $H$ be $\{x, b_1, \ldots, b_k, c_1, \ldots, c_k\}$. Let the edges of $H$ consist of an edge from $x$ to every vertex (including itself), a directed clique on $B = \{b_1, \ldots, b_k\}$ and a directed clique on $C = \{c_1, \ldots, c_k\}$. Let $X$ be the singleton set $\{x\}$. The $H$-colorings of $G$ correspond to the length-$n$ strings satisfying the regular expression $X^*B^* \cup X^*C^*$. Let $A$ be the set of $H$-colorings satisfying the regular expression $X^*B^+$. (A coloring in $A$ starts out with a possibly empty sequence of color-$x$ vertices, then contains a nonempty sequence of vertices with colors from $B$.) Let $M$ be the set of all colorings with at most one color from $B \cup C$. Since $B$ and $C$ are the same size, $\pi(A) \le 1/2$. Furthermore, for $\sigma \in A \backslash M$ and $\tau \in \overline{A} \backslash M$, $P_{\mathrm{Gl}}(\sigma, \tau) = 0$ and $P_{\rightarrow}(\sigma, \tau) = 0$. Claim 2.3 of [12] shows that the mixing time of both of these chains is at least $\pi(A)/8\pi(M)$. Now

$$\pi(A) = \frac{|A|}{|\Omega|} = \frac{|A|}{2|A| + 1} \ge \frac{|A|}{3|A|} = \frac{1}{3}.$$

Also,

$$\pi(M) = \frac{|M|}{|\Omega|} = \frac{2k + 1}{|\Omega|} \le \frac{2k + 1}{k^n},$$

which completes the proof.

**10. Comparisons of scan and Glauber for general graphs.** From the results obtained so far, it seems as if one sweep of systematic scan is equivalent to a linear number of Glauber updates. In the majority of cases examined (Sections 4–7), we have obtained tight asymptotic bounds, and we know the equivalence is exact. Where we don't have tight bounds (Section 8), at least our results are consistent with this supposed equivalence. It is natural to wonder whether a result can be framed that relates scan and Glauber in a more general setting, where the graph $G$ is arbitrary.

In this section we use the comparison method of Diaconis and Saloff-Coste [9] to compare the optimal Poincaré constant $\lambda(\mathcal{M}_{\mathrm{Gl}})$ of Glauber dynamics to the optimal Poincaré constant $\lambda(\mathcal{M}_{\rightarrow})$ of scan. Ideally, we might hope for $\lambda(\mathcal{M}_{\rightarrow}) = \Theta(n\lambda(\mathcal{M}_{\mathrm{Gl}}))$. In fact, the best bounds we can prove lose a factor $n$ in either direction so we have a lower bound for $\lambda(\mathcal{M}_{\rightarrow})$ of $\Omega(\lambda(\mathcal{M}_{\mathrm{Gl}}))$, and an upper bound of $O(n^2\lambda(\mathcal{M}_{\mathrm{Gl}}))$. Moreover, for the lower bound, we need to assume $G$ has bounded degree.

10.1. *Comparing scan to Glauber.*

THEOREM 30. *Suppose $G$ has maximum degree $\Delta$. Let $\mathcal{M}_{\mathrm{Gl}}$ and $\mathcal{M}_{\rightarrow}$ be the Glauber dynamics or systematic scan applied to $H$-colorings of $G$ for a fixed but arbitrary $H$. Then $\lambda(\mathcal{M}_{\mathrm{Gl}}) \le 4q^{\Delta+1}\lambda(\mathcal{M}_{\rightarrow})$.*

PROOF. Suppose $\sigma \to \sigma'$ is a possible Glauber transition, that is, $P_{\mathrm{Gl}}(\sigma, \sigma') > 0$. Let $i$ be the unique vertex satisfying $\sigma_i \ne \sigma'_i$. Say that $\tau \in \Omega$ is *between* $\sigma$ and $\sigma'$ if $\sigma \to \tau$ is a possible scan transition, and additionally: (i) $\tau_i = \sigma'_i$ and (ii) $\tau_j = \sigma_j$ for all $j \sim i$, where "$\sim$" denotes adjacency in $G$. Denote by $B(\sigma, \sigma')$ the set of states between $\sigma$ and $\sigma'$. Consider a scan transition from state $\sigma$, and denote by $\mathcal{E}_i$ the event that, for all $k \in \{i\} \cup \{j : j \sim i\}$, the color proposed by Metropolis($k$) is $\sigma'(k)$. Similarly, consider a reverse scan transition from state $\sigma'$, and denote by $\mathcal{F}_i$ the event that, for all $k \in \{i\} \cup \{j : j \sim i\}$, the color proposed by Metropolis($k$) is $\sigma'(k)$.

The following observations are easy to verify:

- Conditioned on $\mathcal{E}_i$, a scan transition from state $\sigma$ is certain to result in a state $\tau \in B(\sigma, \sigma')$.
- For all $\tau \in B(\sigma, \sigma')$, $P_{\leftarrow}(\tau, \sigma' | \mathcal{F}_i) = P_{\rightarrow}(\sigma, t | \mathcal{E}_i)$.
- $\Pr(\mathcal{E}_i) \ge q^{-(\Delta+1)}$ and $\Pr(\mathcal{F}_i) \ge q^{-(\Delta+1)}$.

It follows from these three observations that

$$\begin{aligned}
&\sum_{\tau \in B(\sigma,\sigma')} \min\{P_{\rightarrow}(\sigma,\tau), P_{\leftarrow}(\tau,\sigma')\} \\
&\quad \ge \sum_{\tau \in B(\sigma,\sigma')} \min\{\Pr(\mathcal{E}_i) P_{\rightarrow}(\sigma,\tau|\mathcal{E}_i), \Pr(\mathcal{F}_i) P_{\leftarrow}(\tau,\sigma'|\mathcal{F}_i)\} \\
&\quad \ge q^{-(\Delta+1)} \sum_{\tau \in B(\sigma,\sigma')} \min\{P_{\rightarrow}(\sigma,\tau|\mathcal{E}_i), P_{\leftarrow}(\tau,\sigma'|\mathcal{F}_i)\} \\
&\quad = q^{-(\Delta+1)} \sum_{\tau \in B(\sigma,\sigma')} P_{\rightarrow}(\sigma,\tau|\mathcal{E}_i) \\
&\quad = q^{-(\Delta+1)}.
\end{aligned}$$

Then for any $f:\Omega\to\mathbb{R}$,

$$\begin{aligned}
\mathcal{EM}_{\mathrm{Gl}}(f,f) &= \frac{1}{2}\sum_{\sigma,\sigma'\in\Omega}\pi(\sigma)P_{\mathrm{Gl}}(\sigma,\sigma')(f(\sigma)-f(\sigma'))^2 \\
&\le \frac{q^{\Delta+1}}{2}\sum_{\sigma,\sigma'\in\Omega}\pi(\sigma)P_{\mathrm{Gl}}(\sigma,\sigma') \\
&\qquad\times\sum_{\tau\in B(\sigma,\sigma')}\min\{P_{\to}(\sigma,\tau),P_{\leftarrow}(\tau,\sigma')\} \\
&\qquad\qquad\times(f(\sigma)-f(\sigma'))^2 \\
&\le q^{\Delta+1}\sum_{\sigma,\sigma'\in\Omega}\pi(\sigma)P_{\mathrm{Gl}}(\sigma,\sigma') \qquad (30)\\
&\qquad\times\sum_{\tau\in B(\sigma,\sigma')}[P_{\to}(\sigma,\tau)(f(\sigma)-f(\tau))^2 \\
&\qquad\qquad+P_{\leftarrow}(\tau,\sigma')(f(\tau)-f(\sigma'))^2] \\
&\le q^{\Delta+1}\sum_{\sigma,\sigma'\in\Omega}\pi(\sigma)P_{\mathrm{Gl}}(\sigma,\sigma') \\
&\qquad\times\sum_{\tau\in B(\sigma,\sigma')}P_{\to}(\sigma,\tau)(f(\sigma)-f(\tau))^2 \\
&\quad+q^{\Delta+1}\sum_{\sigma',\sigma\in\Omega}\pi(\sigma')P_{\mathrm{Gl}}(\sigma',\sigma) \\
&\qquad\times\sum_{\tau\in B(\sigma',\sigma)}P_{\leftarrow}(\tau,\sigma')(f(\tau)-f(\sigma'))^2 \\
&= 2q^{\Delta+1}\sum_{\sigma,\sigma'\in\Omega}\pi(\sigma)P_{\mathrm{Gl}}(\sigma,\sigma')\sum_{\tau\in B(\sigma,\sigma')}P_{\to}(\sigma,\tau)(f(\sigma)-f(\tau))^2 \qquad (31)\\
&= 2q^{\Delta+1}\sum_{\sigma,\tau\in\Omega}\pi(\sigma)P_{\to}(\sigma,\tau)(f(\sigma)-f(\tau))^2\sum_{\sigma':\tau\in B(\sigma,\sigma')}P_{\mathrm{Gl}}(\sigma,\sigma') \\
&\le 2q^{\Delta+1}\sum_{\sigma,\tau\in\Omega}\pi(\sigma)P_{\to}(\sigma,\tau)(f(\sigma)-f(\tau))^2 \qquad (32)\\
&= 4q^{\Delta+1}\mathcal{EM}_{\to}(f,f).
\end{aligned}$$

Inequality (30) applies the fact that $\frac{1}{2}(a-b)^2\le(a-\xi)^2+(\xi-b)^2$ for all $\xi$. Inequality (31) uses the fact that Glauber is time reversible, that is, that $\pi(\sigma)P_{\mathrm{Gl}}(\sigma,\sigma')=\pi(\sigma')P_{\mathrm{Gl}}(\sigma',\sigma)$, for all $\sigma,\sigma'\in\Omega$ and the fact that $B(\sigma,\sigma')=B(\sigma',\sigma)$. Inequality (32) seems crude at first sight, but it is not obvious how to do better: the knowledge of $\tau$ does little to constrain $\sigma'$. $\square$

The inverse of $\lambda(\mathcal{M})$ is closely related to the mixing time of $\mathcal{M}$. Much is known about the precise relationship between these quantities; see, for example, the inequalities in [1, 9, 10, 12, 16, 23, 26]. Some known results only apply when $\mathcal{M}$ is reversible, or when the eigenvalues of its transition matrix $P$ are positive. Our survey paper [14] gives inequalities between Poincaré constants and mixing times in both the general case and the reversible case. We will not repeat the details or trace the development of the ideas here, but we mention a few simple facts that are useful for us. Slightly stronger bounds can be obtained with more effort. Let $\widetilde{\mathcal{M}}_{\rightarrow}$ be the continuization of $\mathcal{M}_{\rightarrow}$ as defined in [2], Chapter 2, page 5. Essentially, this is just $\mathcal{M}_{\rightarrow}$ except that the holding time between discrete transitions is exponential with mean 1. It is a classical result (see, e.g., [23], pages 55 and 63) that the mixing time of $\widetilde{\mathcal{M}}_{\rightarrow}$ is bounded as follows:

$$\mathrm{Mix}_x(\widetilde{\mathcal{M}}_{\rightarrow}, \varepsilon) \le \frac{1}{\lambda(\mathcal{M}_{\rightarrow})}(2\ln(1/\varepsilon) + \ln(1/\pi(x))). \tag{33}$$

Combining (33), Theorem 30 and the upper bound $1/\lambda(\mathcal{M}_{\mathrm{Gl}}) = O(n^4)$ from Section 8.1, we get the following:

THEOREM 31. *Let $H$ be a fixed connected graph. Let $G$ be the $n$-vertex path. Let $\Omega'$ be the state space of $\mathcal{M}_{\mathrm{Gl}}$. It is either the set of all proper $H$-colorings of $G$ (if $H$ is not bipartite) or it is one of the two maximum sets of compatible colorings (if $H$ is bipartite). Consider the Markov chain $\widetilde{\mathcal{M}}_{\rightarrow}$ on the state space $\Omega'$:*

$$\mathrm{Mix}(\widetilde{\mathcal{M}}_{\rightarrow}, \varepsilon) = O(n^5 \ln \varepsilon^{-1}).$$

Let $\mathcal{M}_{\mathrm{Gl}}^{\mathrm{ZZ}}$ be the "lazy" version of Glauber dynamics from page 53 of [23]. In each step, the lazy Markov chain stays where it is with probability 1/2, and otherwise makes the transition specified in the definition of $\mathcal{M}_{\mathrm{Gl}}$. We introduce the lazy chain to keep the eigenvalues positive. See [14] for inequalities which avoid this device. The following inequality from [14] is similar to Proposition 1(ii) of [26]:

$$\frac{1}{\lambda(\mathcal{M}_{\mathrm{Gl}})} = \frac{1}{\lambda(\mathcal{M}_{\mathrm{Gl}}^{\mathrm{ZZ}})} \le \max_x \mathrm{Mix}_x\left(\mathcal{M}_{\mathrm{Gl}}^{\mathrm{ZZ}}, \frac{1}{2e}\right).$$

Combining this with Theorem 30 and with (33), we find that, for bounded-degree graphs $G$, the mixing time of $\widetilde{\mathcal{M}}_{\rightarrow}$ is at most $O(n)$ times the mixing time of $\mathcal{M}_{\mathrm{Gl}}^{\mathrm{ZZ}}$. Perhaps this result can be improved by a factor of $n^2$.

10.2. *Comparing Glauber to scan.*

THEOREM 32. *Suppose $G$ is arbitrary. Let $\mathcal{M}_{\mathrm{Gl}}$ and $\mathcal{M}_{\rightarrow}$ be the Glauber dynamics and systematic scan applied to $H$-colorings of $G$ for a fixed but arbitrary $H$. Then $\lambda(\mathcal{M}_{\rightarrow}) \le n^2 q \lambda(\mathcal{M}_{\mathrm{Gl}})$.*

PROOF. Let $\sigma, \sigma' \in \Omega$ be a pair of states for which $P_{\rightarrow}(\sigma, \sigma') > 0$. There is a natural canonical path $\gamma_{\sigma,\sigma'} = (\sigma = \tau^0 \rightarrow \tau^1 \rightarrow \cdots \rightarrow \tau^n = \sigma')$ from $\sigma$ to $\sigma'$ using Glauber transitions, in which $\tau^{i-1}$ differs from $\tau^i$ (if at all) only at vertex $i$. According to [9], Theorem 2.1, the quantity we need to bound is

$$\begin{aligned} A &= \frac{1}{\pi(\tau) P_{\mathrm{Gl}}(\tau, \tau')} \sum_{\sigma, \sigma' : (\tau, \tau') \in \gamma_{\sigma,\sigma'}} \pi(\sigma) P_{\rightarrow}(\sigma, \sigma') |\gamma_{\sigma,\sigma'}| \\ &= n^2 q \sum_{\sigma, \sigma' : (\tau, \tau') \in \gamma_{\sigma,\sigma'}} P_{\rightarrow}(\sigma, \sigma'), \end{aligned}$$

where we have used the facts that $\pi$ is uniform, $|\gamma_{\sigma,\sigma'}| = n$, and $P_{\mathrm{Gl}}(\tau, \tau') = 1/nq$. (Diaconis and Saloff-Coste state their theorem for time-reversible MCs, but their proof does not use time reversibility.) We shall demonstrate that $A \le n^2 q$, from which it follows that $\lambda(\mathcal{M}_{\rightarrow}) \le n^2 q \lambda(\mathcal{M}_{\mathrm{Gl}})$.

Regard $\tau$ and $\tau'$ as fixed, and suppose $\tau$ and $\tau'$ differ at vertex $i$. Denote by $\mathcal{E}_1^{\sigma}$ the event that the sequence

$$\mathrm{Metropolis}(1), \mathrm{Metropolis}(2), \ldots, \mathrm{Metropolis}(i-1)$$

takes $\sigma$ to $\tau$, and by $\mathcal{E}_2^{\sigma'}$ the event that

$$\mathrm{Metropolis}(i+1), \mathrm{Metropolis}(i+2), \ldots, \mathrm{Metropolis}(n)$$

takes $\tau'$ to $\sigma'$. Then $P_{\rightarrow}(\sigma, \sigma') \le \Pr(\mathcal{E}_1^{\sigma} \wedge \mathcal{E}_2^{\sigma'}) = \Pr(\mathcal{E}_1^{\sigma}) \Pr(\mathcal{E}_2^{\sigma'})$, and

$$\sum_{\sigma, \sigma' : (\tau, \tau') \in \gamma_{\sigma,\sigma'}} P_{\rightarrow}(\sigma, \sigma') \le \sum_{\sigma} \Pr(\mathcal{E}_1^{\sigma}) \sum_{\sigma'} \Pr(\mathcal{E}_2^{\sigma'}).$$

The second sum above is clearly bounded by 1, since the events $\{\mathcal{E}_2^{\sigma'} : \sigma' \in \Omega\}$ are disjoint. In fact, the first sum is also bounded by 1, since $\Pr(\mathcal{E}_1^{\sigma})$ is equal to the probability that the sequence

$$\mathrm{Metropolis}(i-1), \mathrm{Metropolis}(i-2), \ldots, \mathrm{Metropolis}(1)$$

takes $\tau$ to $\sigma$. So the terms of the first sum may also be viewed as probabilities of disjoint events. Thus, $A \le n^2 q$, as claimed. $\square$

Combining Theorem 32 with inequalities of Diaconis and Stroock [10] and Sinclair [26], we get

$$\mathrm{Mix}_x(\mathcal{M}_{\mathrm{Gl}}^{\mathrm{ZZ}}, \varepsilon) \le n^2 q \frac{1}{\lambda(\mathcal{M}_{\rightarrow})} \ln \frac{1}{\varepsilon \pi(x)}. \tag{34}$$

This can be combined with the upper bound

$$\frac{1}{\lambda(\mathcal{M}_{\rightarrow})} \leq \frac{2(\max_x \mathrm{Mix}_x(\mathcal{M}_{\rightarrow}, 1/e))^2}{(1/2 - 1/e)^2}. \tag{35}$$

The square of the mixing time in (35) is necessary in the general nonreversible case (see [14]), though of course better inequalities might apply to the particular chain $\mathcal{M}_{\rightarrow}$. Combining (34) and (35), we get a weak inequality which shows that the mixing time of (lazy) Glauber dynamics is at most $O(n^3)$ times the square of the mixing time of systematic scan.

Note that the proofs of Theorems 30 and 32 are actually about Dirichlet forms rather than about Poincaré constants, so the same inequalities apply to the log-Sobolev constant.

**Acknowledgment.** The authors thank Persi Diaconis for comments on a draft of this article.

M. Dyer
School of Computing
University of Leeds
Leeds LS2 9JT
UK
url: www.comp.leeds.ac.uk/~dyer/

L. Goldberg
Department of Computer Science
University of Warwick
Coventry CV4 7AL
UK
url: www.dcs.warwick.ac.uk/~leslie/

M. Jerrum
School of Informatics
University of Edinburgh
JCMB, The King's Buildings
Edinburgh EH9 3JZ
UK
e-mail: mrj@inf.ed.ac.uk
url: homepages.inf.ed.ac.uk/mrj/